\documentclass[12pt,francais, reqno]{smfart}

\usepackage{amsmath}
\usepackage{amssymb}
\usepackage{amsthm}
\usepackage{a4wide}
\usepackage[frenchb,english]{babel}
\usepackage{mathrsfs}
\usepackage{graphics}
\usepackage{latexsym}

\numberwithin{equation}{section}
\font\xrm=wncyr10
\def\shu{\hbox{\xrm sh} }

\newcommand{\lm}{\left}
\newcommand{\rt}{\right}

\newcommand{\un}{\underline}
\newcommand{\dd}{\textup{d}}
\newcommand{\ba}{\textup{B}}

\newcommand{\la}{\textup{La}}
\newcommand{\li}{\textup{Li}}

\def\al{\alpha}

\def\Li{\operatorname{Li}}
\def\dis{\displaystyle}
\def\dd{\textup{d}}

\def\di{\displaystyle}
\def\bs{\underline{s}}

\def\wzeta{\overline{\zeta}}

\title[]{S\'eries hyperg\'eom\'etriques multiples et polyz\^etas}
\author[]{J. Cresson, S. Fischler et T. Rivoal}

\date{30 juin 2006}

\newtheorem{theo}{Th\'eor\`eme}
\newtheorem{lem}{Lemme}
\newtheorem{prop}{Proposition}

\newtheorem{conj}{Conjecture}
\newtheorem{remark}{Remarque}

\newcommand{\gdo}{\mathcal{O}}

\newcommand{\eps}{\varepsilon}

\newcommand{\C}{\mathbb C}

\newcommand{\indso}[2]{{\tiny {\begin{array}{c} #1 \\ #2 \end{array}}}}
\newcommand{\indtr}[3]{{\tiny {\begin{array}{c} #1 \\ #2 \\ #3 \end{array}}}}
\newcommand{\Ic}{I^{{\rm c}}}
\newcommand{\unp}{\{1, \ldots, p\}}
\newcommand{\moins}{\setminus}
\newcommand{\sich}{\hat s_i}
\newcommand{\stch}{\hat s_t}
\newcommand{\aich}{\hat A_i}
\newcommand{\Ciandco}{C\left[\,{\tiny \begin{matrix} I \\ (s_i) \\ (j_i) \\ (\sich) \end{matrix}} \,\right]}
\newcommand{\diandco}{\partial\left[\,{\tiny \begin{matrix} I \\ (s_i) \\  (\sich) \end{matrix}} \,\right]}

\newcommand{\qqq}{\varpi}
\newcommand{\cale}{{\mathscr{E}}}
\newcommand{\Rtilde}{\widetilde R}
\newcommand{\Qtilde}{\widetilde Q}
\newcommand{\Pti}{\widetilde P}

\begin{document}
\maketitle

\setcounter{tocdepth}{2}
\baselineskip 6mm

\begin{abstract}
Nous d\'ecrivons un algorithme th\'eorique et effectif permettant  de d\'emontrer que
des s\'eries et int\'egrales hyperg\'eom\'e\-tri\-ques multiples relativement g\'en\'erales se
d\'ecomposent en combinaisons lin\'eaires \`a coefficients rationnels de
polyz\^etas.
\end{abstract}

\begin{altabstract}
We describe a theoretical and effective algorithm which enables
us to prove that rather general hypergeometric series and integrals can be decomposed
as linear combinations of multiple zeta values, with rational coefficients.
\end{altabstract}

\tableofcontents

\section{Introduction}\label{sec:1}
Une g\'en\'eralisation
de la fonction z\^eta de Riemann $\zeta(s)$ est
donn\'ee par les s\'eries {\em polyz\^etas}, d\'efinies pour tout entier
$p\ge 1$
et tout $p$-uplet $\underline{s}=(s_1, s_2, \dots, s_p)$
d'entiers $\ge 1$,
avec $s_1\ge 2$, par
$$
\zeta(s_1, s_2, \ldots, s_p)=
\sum_{k_1> k_2>\cdots > k_p\ge 1}
\frac{1}{k_1^{s_1}k_2^{s_2}\cdots k_p^{s_p}}.
$$
Les entiers $p$ et
$s_1+s_2+\cdots+s_p$ sont respectivement la profondeur et
le poids de $\zeta(s_1, s_2, \ldots, s_p)$. Pour diverses
raisons, il est plus simple de consid\'erer que la sommation
est faite sur  $k_1\ge k_2\ge \cdots\ge k_p \ge 1$ : nous
noterons $\wzeta(s_1, s_2, \ldots, s_p)$ les s\'eries ainsi
obtenues.  Il est \`a noter que les deux s\'eries
convergent plus g\'en\'eralement
pour des exposants complexes v\'erifiant
$\sum _{j=1} ^{r}\Re(s_{j})> r$ pour tout $r\in\{1,\dots,p\}$, ce qui
autorise \`a avoir des exposants entiers n\'egatifs par exemple.

Les polyz\^etas interviendront dans cet article par l'interm\'ediaire des
fonctions polylogarithmes multiples, d\'efinies par
$$
\Li_{s_1, s_2 \ldots, s_p}(z_1, z_2,\ldots, z_p)=
\sum_{ k_1> k_2 >\cdots > k_p\ge 1}
\frac{z_1^{k_1}z_2^{k_2}\cdots z_p^{k_p}}{k_1^{s_1}k_2^{s_2}
\cdots k_p^{s_p}}
$$
pour $\vert z_1\vert\le 1, \ldots,
\vert z_p\vert \le 1$.
On obtiendra en fait
les r\'esultats  pour
les
polylogarithmes multiples larges, d\'efinis par
$$
\la_{s_1, s_2\ldots, s_p}(z_1, z_2\ldots, z_p)=
\sum_{ k_1 \ge k_2\ge \cdots \ge k_p\ge 1}
\frac{z_1^{k_1}z_2^{k_2}\cdots z_p^{k_p}}{k_1^{s_1}k_2^{s_2}
\cdots k_p^{s_p}}.
$$
Lorsque $p=1$, les deux variantes co{\"\i}ncident avec les polylogarithmes
usuels et si $z_1=z_2=\cdots =z_p=1$ et $s_1 \geq 2$, on a
$\Li_{s_1, s_2\ldots, s_p}(1, 1, \ldots, 1)=\zeta(s_1,s_2\ldots, s_p)$ et
$\la_{s_1, s_2\ldots, s_p}(1,1, \ldots, 1)=\wzeta(s_1,s_2\ldots, s_p).$
Un th\'eor\`eme d'Ulanski{\u\i}~\cite{ulanskii}
permet de passer lin\'eairement d'un type de s\'erie
\`a l'autre ; en vue d'applications diophantiennes,
on ne perd donc rien \`a consid\'erer une variante plut\^ot qu'une autre.

Remarquons d\`es \`a pr\'esent que les fonctions polylogarithmes multiples peuvent \^etre
d\'efinies pour des exposants $s_i$ complexes, \`a condition de supposer en plus que
$\vert z_1\vert <1$ pour des raisons de convergence. En particulier, nous utiliserons
ces fonctions avec des $s_j \in \mathbb{Z}$ : par d\'efinition,
le poids d'une telle fonction est alors
$\sum_{j=1}^p \max(s_j,0).$

On voit naturellement appara{\^\i}tre les polyz\^etas lorsque, par exemple,
on consid{\`e}re les produits des valeurs de la fonction
z\^eta :
on a  $\zeta(n)\zeta(m)=\zeta(n+m)+\zeta(n,m)+\zeta(m,n)$,
ce qui permet en
quelque sorte  de {\og lin{\'e}ariser\fg} ces produits.
En dehors de quelques identit\'es telles que $\zeta(2,1)=\zeta(3)$
 (due \`a Euler),
la nature arithm\'etique de ces s\'eries est aussi peu connue que
celle des nombres $\zeta(s)$.
Cependant, l'ensemble des nombres $\zeta(\underline s)$ poss\`ede
une tr\`es riche
structure alg\'ebrique assez bien comprise, au moins
conjecturalement (voir~\cite{miw}).
Par exemple, on peut s'int{\'e}resser
 aux $\mathbb{Q}$-sous-espaces vectoriels $\mathcal{Z}_p$ de
 $\mathbb{R}$, engendr{\'e}s par les $2^{p-2}$ polyz{\^e}tas de
poids $p\ge 2$ :  $\mathcal{Z}_2=\mathbb{Q}\zeta(2)$,
$\mathcal{Z}_3=\mathbb{Q}\zeta(3)+\mathbb{Q}\zeta(2,1)$,
$\mathcal{Z}_4=\mathbb{Q}\zeta(4)+\mathbb{Q}\zeta(3,1)+
\mathbb{Q}\zeta(2,2)+
\mathbb{Q}\zeta(2,1,1)$,
etc. Posons $v_p=\textup{dim}_{\mathbb{Q}}(\mathcal{Z}_p)$.
On a alors la

\begin{conj}\label{conj:KZ}
$(i)$
Pour tout entier $p\ge 2$, on a $v_p=c_p$, o{\`u} l'entier
$c_p$ est d{\'e}fini par
la r{\'e}currence de type Fibonacci
$c_{p+3}=c_{p+1}+c_{p}$,
avec $c_0=1$, $c_1=0$ et $c_2=1$.

 $(ii)$ Les $\mathbb{Q}$-espaces vectoriels $\mathbb{Q}$ et
$\mathcal{Z}_p$ ($p\ge 2)$, sont en somme
directe.
\end{conj}
La suite $(v_p)_{p\ge 2}$ devrait donc cro{\^\i}tre comme
$\al^p$ (o{\`u} $\al\approx
1,3247$ est racine du
polyn{\^o}me $X^3-X-1$), ce qui est  bien plus petit que
$2^{p-2}$. Il y a
donc conjecturalement beaucoup de relations lin{\'e}aires entre les polyz{\^e}tas
de m{\^e}me poids
et aucune en poids diff{\'e}rents :
dans cette direction, un
th{\'e}or{\`e}me de Goncharov~\cite{goncharov1} et Terasoma~\cite{terasoma}
affirme que l'on a  $v_p\le c_p$ pour tout entier $p\ge 2$. Il
reste donc {\`a} montrer l'in{\'e}galit{\'e} inverse pour montrer (i) mais aucune
minoration non
triviale de $v_p$ n'est connue {\`a} ce jour : si l'on montre
facilement que $v_2=v_3=v_4=1$,
on est bloqu{\'e} d{\`e}s l'{\'e}galit{\'e} $v_5=2$, qui
est {\'e}quivalente
{\`a} l'irrationalit{\'e} toujours inconnue de
$\zeta(5)/(\zeta(3)\zeta(2))$.
Plus g\'en\'eralement, un des int\'er\^ets de la
conjecture~\ref{conj:KZ} est d'impliquer la
suivante.

\begin{conj}\label{conj:zeta alg inde}
Les nombres $\,\pi, \zeta(3), \zeta(5), \zeta(7), \zeta(9),$ etc,
sont alg{\'e}brique\-ment ind{\'e}pen\-dants sur $\mathbb{Q}$.
\end{conj}

Cette conjecture semble actuellement totalement hors de
port\'ee. Un certain nombre de r\'esultats diophantiens ont
n\'eanmoins
\'et\'e  obtenus en profondeur 1, c'est-\`a-dire dans le cas
de la
fonction z\^eta de Riemann :
\begin{itemize}
\item[(i)] Le nombre $\zeta(3)$ est irrationnel (Ap\'ery~\cite{ap}) ;
\item[(ii)] La dimension de l'espace vectoriel engendr\'e sur $\mathbb{Q}$ par 1,
$\zeta(3)$, $\zeta(5), \ldots, \zeta(A)$ (avec $A$ impair) cro{\^\i}t au moins comme $\log(A)$
(\cite{br, ri}) ;
\item[(iii)] Au moins un des quatre nombres $\zeta(5), \zeta(7), \zeta(9), \zeta(11)$ est irrationnel
(Zudilin~\cite{zud}).
\end{itemize}
Ces r\'esultats peuvent \^etre obtenus par l'\'etude de
certaines s\'eries de la forme
\footnote{du moins, dans le cas de (i) et (ii) ; le point (iii)
n\'ecessite une id\'ee {\em a priori} diff\'erente (s\'erie {\og d\'eriv\'ee\fg})
mais on peut l'int\'egrer dans le cadre fourni par~\eqref{eq:serie simple generale}.
Voir un peu plus loin dans cette Introduction pour plus de d\'etails.}
\begin{equation}\label{eq:serie simple generale}
\sum_{k=1}^{\infty} \frac{P(k)}{k^A(k+1)^A\cdots (k+n)^A} \,z^{-k}
\end{equation}
avec $P(X)\in\mathbb{Q}[X]$, $n\ge0$, $A\ge 1$ et
$\vert z \vert\ge 1$ (le choix de $1/z$ plut\^ot que $z$
est purement technique) :
nous rappelons
sommairement au paragraphe~\ref{ssec:serieunevar} comment on
utilise ces s\'eries pour les d\'emontrer,
en exploitant le fait que, g\'en\'eriquement, elles s'expriment
aussi comme combinaisons  lin\'eaires des valeurs de z\^eta aux
entiers lorsque $z=1$. Les divers choix de $P$ conduisent \`a
des {\em s\'eries hyperg\'eom\'etriques g\'en\'eralis\'ees} : voir
les ouvrages~\cite{ba, sl} pour les d\'efinitions, qui ne sont pas
essentielles ici.

Notre but est de poser les bases d'une g\'en\'eralisation de
cette m\'ethode
hyperg\'eom\'etri\-que en profondeur quelconque en consid\'erant
{\em a priori} des s\'eries multiples de la forme
\begin{equation}
\label{eq:serie mult generale}
\sum_{k_1\ge \cdots \ge k_p\ge 1}
\frac{P(k_1, \ldots, k_p)}{(k_1)_{n_1+1}^{A_1}
\cdots (k_p)_{n_p+1}^{A_p}}\,z_1^{-k_1}
\cdots z_p^{-k_p},
\end{equation}
avec $P(X_1,\ldots, X_p)\in\mathbb{Q}[X_1,\ldots, X_p]$, des
entiers $A_j\ge 1$ et $n_j\ge 0$ et
$\vert z_1\vert\ge 1, \ldots,
\vert z_p\vert \ge 1$,  ceci dans l'espoir qu'elles
s'expriment comme combinaisons lin\'eaires de polyz\^etas
int\'eres\-sants lorsque
$z_1= \cdots =z_p=1$. (Pour raccourcir les
expressions, on a utilis\'e le symbole de Pochhammer
$(\al)_m=\al(\al+1)\cdots (\al+m-1)$.) On pourrait imaginer g\'en\'eraliser encore
\eqref{eq:serie mult generale} en rempla\c cant, au d\'enominateur, chaque facteur $(k_i)_{n_i+1}^{A_i}$
par $(k_i+r_i)_{n_i+1}^{A_i}$. Cela peut \^etre utile (et nos m\'ethodes le permettent) si des bornes explicites apparaissent en fonction des $n_i$, mais pour des r\'esultats qualitatifs c'est inutile car on peut s'y ramener, en rempla\c cant $n_i$ par $n_i + r_i$ et en multipliant le num\'erateur par $(k_i)_{r_i} ^{A_i}$.

Les  s\'eries de la forme \eqref{eq:serie mult generale}  apparaissent naturellement dans la litt\'erature.
Par exemple, Sorokin~\cite{so1} a d\'eduit l'irrationalit\'e
de $\zeta(3)$ d'un
r\'esultat que l'on peut \'ecrire ainsi (voir \S \ref{sec:prop:intesoro=seriemult})~: pour tout entier
$n\ge 0$, on a
\begin{equation}
\label{eq:sorokin1}
n! \sum_{k_1\ge k_2\ge 1}
 \frac{(k_2-n)_n(k_1-k_2+1)_n }{(k_1)_{n+1}^2(k_2)_{n+1}} =
2a_n\zeta(2,1)-b_n,
\end{equation}
o\`u $a_n$ et $b_n$ sont les c\'el\`ebres nombres rationnels
utilis\'es par Ap\'ery~\cite{ap} dans sa
preuve originelle de l'irrationalit\'e de $\zeta(3)$. La
m\'ethode de Sorokin n'utilise pas directement la s\'erie~\eqref{eq:sorokin1}
 mais
consiste \`a  r\'esoudre un subtil probl\`eme d'approximation
de Pad\'e, qu'il n'est malheureusement
pas facile de g\'en\'eraliser \`a d'autres
situations.
Nous nous affranchissons de
l'approximation de Pad\'e pour esp\'erer profiter, en
profondeur sup\'erieure, de la grande souplesse de la m\'ethode
hyperg\'eom\'etrique en profondeur 1. Il
est int\'eressant de noter que la s\'erie double en~\eqref{eq:sorokin1}
est un exemple de {\em s\'erie hyperg\'eom\'etrique de Kamp\'e de F\'eriet}
(voir~\cite[p. 27]{srivastava}), comme on le voit apr\`es
quelques transformations triviales du sommande.
Par un l\'eger abus de langage,
nous appelons s\'erie hyperg\'eom\'etrique multiple une expression
de la forme~\eqref{eq:serie mult generale} bien que, en g\'en\'eral, il ne s'agisse seulement que de
combinaisons lin\'eaires rationnelles
de telles s\'eries.

Un ingr\'edient, fr\'equemment utilis\'e avec des s\'eries simples, consiste \`a
d\'eriver la fraction rationnelle en $k$ dans~la
s\'erie~\eqref{eq:serie simple generale}, avant de sommer ; par exemple, une double d\'erivation sert \`a d\'emontrer le r\'esultat de
Zudilin~\cite{zud} rappel\'e apr\`es la conjecture~\ref{conj:zeta alg inde}.
Cette astuce, appliqu\'ee plusieurs fois, permet de faire dispara\^{\i}tre $\zeta(s)$ de la forme
lin\'eaire obtenue, pour de petites valeurs de $s$. On peut imaginer l'utiliser pour des sommes multiples, m\^eme
si on n'a aucun r\'esultat connu de disparition de polyz\^etas dans ce cadre. Il est clair qu'en d\'erivant
une fraction rationnelle de la forme
$P(X_1, \ldots, X_p)/\big((X_1)_{n+1}^{A_1} \ldots (X_p)_{n+1}^{A_p} \big)$ par rapport \`a l'une des variables $X_j$,
on obtient une fraction rationnelle de la m\^eme forme (avec $A_j$ remplac\'e par $A_j+1$) : cette
remarque montre que l'on ne perd rien \`a consid\'erer des s\'eries de la forme~\eqref{eq:serie mult generale}.

En profondeur $p\ge2$, l'\'etude des s\'eries multiples du type
de~\eqref{eq:serie mult generale} se d\'ecompose
en plusieurs \'etapes et, malheureusement, la
premi\`ere difficult\'e se pr\'esente d\`es la premi\`ere
\'etape, qui est pourtant triviale en profondeur 1. Nous
mettons ceci en \'evidence sur l'exemple de la profondeur 2
au paragraphe~\ref{ssec: algo prof 2} : la g\'en\'eralisation en profondeur
quelconque n\'ecessite la production d'un algorithme r\'ecursif
(permettant de d\'eduire le cas de la profondeur $p$ du cas de
la profondeur $p-1$) que l'on d\'ecrit au paragraphe~\ref{sec:2}.
Informellement, on obtient alors le r\'esultat suivant.
\begin{theo}\label{theo:informel}
Supposons que l'on ait
$\vert z_1 \vert >1$ et
$\vert z_j \vert \ge  1$ pour tout $j=2,\dots ,p$.
Alors, toute s\'erie de la
forme~\eqref{eq:serie mult generale} s'\'ecrit comme une
combinaison lin\'eaire \`a coefficients polyn\^omes de Laurent
dans $\mathbb{Q}[z_1^{\pm 1}, \ldots, z_p^{\pm 1}]$ en les
polylogarithmes multiples
$\la_{s_1, \ldots, s_q}(1/\widehat z_1,\ldots, 1/\widehat z_q)$
o\`u $0\le q\le p$, $\sum_{j=1} ^q \max(s_j,0) \le \sum_{j=1} ^p A_j$ et o\`u
les $\widehat z_1,\ldots, \widehat z_q$ sont certains
produits des $z_1, \ldots, z_p.$
\end{theo}
\begin{remark} $(1)$ Bien que peu surprenant en apparence, ce r\'esultat
est, comme on le verra, loin d'\^etre facile \`a d\'emontrer. Ici, les entiers $s_1$, \ldots, $s_q$ peuvent \^etre de signe quelconque et, comme dans toute la suite,
on doit entendre un polylogarithme de profondeur $0$
comme \'etant la fonction identiquement \'egale \`a $1$.

$(2)$ Certains des $s_j$ peuvent \^etre n\'egatifs ou nuls : cela ne peut \^etre le cas que si
l'un des degr\'es en l'une des variables $X_j$ de la fraction
$P(X_1,  \ldots, X_p)/
\big((X_1)_{n_1+1}^{A_1} \cdots (X_p)_{n_p+1}^{A_p}\big)$ est
positif, c'est-\`a-dire lorsque $\textup{deg}_{X_j}(P)\ge A_j(n_j+1)$.

$(3)$ On peut raffiner ce th\'eor\`eme~: voir
le Th\'eor\`eme~\ref{theo:raffinement} au paragraphe~\ref{sec:raffinement algo}.
Il en r\'esulte par exemple que les
polyn\^omes de Laurent sont en fait toujours dans
$\mathbb{Q}[z_1, z_2^{\pm 1}, \ldots, z_p^{\pm 1}]$.
\end{remark}

Une deuxi\`eme difficult\'e provient du fait que certains
polylogarithmes multiples peuvent avoir un ou des exposants $s_j
\le 0$, ce qui n\'ecessite un traitement \`a part. On obtient
le r\'esultat dit de {\it non-enrichissement} suivant  (voir le
paragraphe~\ref{sec:nonenrichissement}): {\em
Lorsque les modules de
$z_{1}, \ldots, z_{p}$ sont tous diff\'erents de $1$, tout  polylogarithme
multiple $\la_{s_1, \ldots, s_p}(\underline z)$, de profondeur $p$
et ayant certains exposants $\le 0$, est
 une combinaison lin\'eaire en des polylogarithmes multiples d'indices
$\ge 1$ (en des produits des $z_j$) de poids $\le \sum_{j=1}^p \max(s_j,0)$,
dont les  coefficients sont des polyn\^omes \`a coefficients rationnels
en les
$\di \big((1-z_{j_1}\cdots z_{j_m})^{-1}\big)_{1\le j_1< \cdots < j_m\le p,\, m\ge 1}$ et
les $(z_j^{\pm 1})_{1\le j\le p}$.}

En combinant ce r\'esultat et le th\'eor\`eme \ref{theo:informel} on obtient l'\'enonc\'e suivant (qui a \'et\'e obtenu ind\'ependamment, dans le cas particulier $z_1 = \ldots = z_p$, par Zlobin \cite{ZlobinZametki2005}) :

\begin{theo}
\label{theofinalenz}
Supposons que pour tout $j=1,\dots ,p$, on ait $ \vert z_j \vert > 1$. Alors, toute s\'erie
de la forme~\eqref{eq:serie mult generale} s'\'ecrit comme une
combinaison lin\'eaire \`a coefficients polyn\^omes \`a coefficients rationnels en les
$\di \big((1-z_{j_1}\cdots z_{j_m})^{-1}\big)_{1\le j_1< \cdots < j_m\le p,\, m\ge 1}$ et
les $(z_j^{\pm 1})_{1\le j\le p}$ de polylogarithmes multiples
$\la_{s_1, \ldots, s_q}(1/\widehat z_1,\ldots, 1/\widehat z_q)$
o\`u $0\le q\le p$, $s_i \geq 1$ pour $i=1,\dots ,q$,
$\sum_{j=1} ^q s_j \le \sum_{j=1} ^p A_j$ et o\`u
les $\widehat z_1,\ldots, \widehat z_q$ sont certains
produits des $z_1, \ldots, z_p.$
\end{theo}

L'analogue des
th\'eor\`emes~\ref{theo:informel} et~\ref{theofinalenz}
lorsque $z_1=\dots =z_p =1$ s'\'enonce
comme suit. Une version plus pr\'ecise (le th\'eor\`eme~\ref{theoreg})
sera d\'emontr\'ee au paragraphe~\ref{ssec:preuvetheoreg}; la d\'emonstration n\'ecessite d'utiliser la r\'egularisation des polyz\^etas divergents.
Ce th\'eor\`eme est celui que nous avons
impl\'ement\'e dans~\cite{algo}, ce qui nous a
permis d'avoir l'id\'ee du th\'eor\`eme \ref{conj:anti-sym} ci-dessous et
d'observer d'autres exemples de s\'eries qui font appara\^\i tre seulement certains des
polyz\^etas attendus \cite{fialgo}.

\begin{theo} \label{theocv}
Toute s\'erie convergente  de la forme (\ref{eq:serie mult generale}) s'\'ecrit
lorsque $z_1=\dots =z_p =1$ comme une combinaison lin\'eaire \`a
coefficients rationnels en les polyz\^etas
$\zeta (s_1 ,\dots ,s_q )$ o\`u $0\leq q\leq p$, $s_1 \geq 2$,
$s_i \geq 1$ pour  $i=1,\dots ,q$ et
$\sum_{j=1}^q s_j \leq \sum_{j=1}^p A_j$.
\end{theo}

Notre algorithme donne diverses pr\'ecisions sur les
th\'eor\`emes~\ref{theo:informel}, \ref{theofinalenz} et \ref{theocv} (d\'enominateurs des
coefficients, degr\'e des polyn\^omes dans le cas des s\'eries les
plus simples, dites briques). De plus, il se pr\^ete (pour tous $z_1$, \ldots, $z_p$) \`a une
impl\'ementation informatique  que nous avons effectu\'ee  \cite{algo} (lorsque $z_1=\dots =z_p =1$)  \`a
l'aide du programme GP/Pari : cela  nous a permis de tester de
nombreuses s\'eries et d'obtenir des r\'esultats tels que
\begin{multline}\label{eq:exemplehorrible}
\sum_{k_1\ge k_2\ge 1} \frac{5k_2^2-k_1^2-4k_1k_2-3k_1+7k_2}{(k_1)_{3}^4\;(k_2+1)_{4}^3}
\\= -\frac{153060027667}{1289945088} + \frac{832127737}{17915904}\, \zeta(2)
+ \frac{33349589}{2985984} \,\zeta(3) + \frac{10561397}{2985984} \,\zeta(4)
\\
+ \frac{117277}{10368} \,\zeta(5)
+ \frac{1475}{1728} \,\zeta(6) + \frac{757}{432} \,\zeta(7)
+ \frac{6125}{1728} \,\zeta(2,2)
\\
+ \frac{245}{24}\,\zeta(2,3) + \frac{35}{32} \,\zeta(3,2)
+ \frac16 \,\zeta(3,3) + \frac{595}{864}\,\zeta(4,2) + \frac74 \,\zeta(4,3).
\\
\end{multline}
Ce r\'esultat pourrait \'eventuellement \^etre un peu simplifi\'e en utilisant les
relations lin\'eaires connues entre polyz\^etas.

Une fois cette \'etape franchie, une troisi\`eme
difficult\'e provient de la profusion de polyz\^etas qui semblent
appara{\^i}tre spontan\'ement dans des exemples {\og au hasard\fg}
comme~\eqref{eq:exemplehorrible}. Nous avons donc \'et\'e conduits \`a
rechercher une classe
de polyn\^omes $P(X_1,\ldots, X_p)$ tels que,
{\em a priori}, seulement certains polyz\^etas
int\'eressants ont un coefficient non-nul \`a la sortie de
l'algorithme. Par {\og int\'eressants \fg}, nous entendons des
polyz\^etas qui ne sont pas trivialement des puissances de $\pi$, qui
parasitent les applications diophantiennes en les rendant
triviales.~\footnote{Par exemple, la minoration de la dimension de
l'espace des nombres $\zeta(2n+1)$ devient
sans int\'er\^et
lorsque l'on
rajoute les nombres $\zeta(2n)$ : la transcendance de $\pi$ implique leur ind\'ependance lin\'eaire
sur $\mathbb{Q}$ et donc une minoration de dimension de l'ordre de $A/2$ au lieu de $\log(A).$
}
Voici quelques exemples de s\'eries qui ne font pas appara\^itre $\pi$ :
\begin{multline*}
\sum_{k_1\ge k_2\ge 1} (k_1+1)(k_2+1)\frac{(k_1-k_2-1)_3(k_1+k_2+1)_3(k_1-1)_5(k_2-1)_5}
{(k_1)_3^5\;(k_2)_3^5}
\\
= \frac{27875}{8192}-\frac{2847}{1024}\,\zeta(3) -\frac{15}{32}\,\zeta(5)+\frac{27}{64}\,\zeta(7),
\end{multline*}
\begin{multline*}
\sum_{k_1\ge k_2\ge 1} \big(k_1+\frac 12\big)\big(k_2+\frac 12\big)\frac{(k_1-k_2-1)_3(k_1+k_2)_3(k_1-1)_4(k_2-1)_4}
{(k_1)_3^7\;(k_2)_3^7}
\\
= -1156 +891\,\zeta(3)+ \frac{189}2 \,\zeta(5) + 78 \big(\zeta(5,3) -\zeta(3,5)\big),
\end{multline*}
\begin{multline*}
\sum_{k_1\ge k_2\ge 1} \frac{(k_1-k_2)(k_1+k_2+4)(k_1-2)_9(k_2-2)_9}
{(k_1)_5^{4}\;(k_2)_5^{4}}
\\
= -\frac{642739948033}{41278242816}+\frac{10214719}{995328}\,\zeta(3)
+ \frac{57497}{18432}\,\zeta(5),
\end{multline*}
\begin{multline*}
\sum_{k_1\ge k_2\ge 3\ge 1}\big(k_1+\frac12\big)\big(k_2+ \frac12\big)\big(k_3+\frac12\big)
\\
\times \frac{(k_1-k_2)(k_2-k_3)(k_1-k_3)
(k_1+k_2+1)(k_1+k_3+1)(k_2+k_3+1)}{(k_1)_2^4\;(k_2)_2^4\;(k_3)_2^4}
\\
= -\frac{1}{4} - \zeta(3) + \frac14 \,\zeta(5) + \zeta(3)^2 -\frac14 \,\zeta(7).
\end{multline*}
Nous avons propos\'e dans \cite{CFRsym} une g\'en\'eralisation en  profondeur quelconque,
des s\'eries {\em very-well-poised}~\footnote{Voir le paragraphe~\ref{ssec:serieunevar} pour
l'origine de cette d\'enomination.} (ou  {\em tr\`es bien \'equilibr\'ees}) introduites en profondeur 1 ; elle
explique les quatre exemples ci-dessus.
 Il s'agit du r\'esultat suivant, qui est d\'emontr\'e sous une forme plus pr\'ecise dans \cite{CFRsym}.

\begin{theo}\label{conj:anti-sym} Fixons trois entiers $A\ge 2$, $n \ge 0$ et $p \ge 1 $,
ainsi que  $P(X_1, \ldots, X_p) \in\mathbb{Q}[X_1, \ldots, X_p]$ un polyn\^ome tel que :
\begin{equation*}
P(X_{\sigma(1)}, X_{\sigma(2)},\ldots, X_{\sigma(p)})
= \varepsilon(\sigma) P(X_1, X_2, \ldots, X_p)
\end{equation*}
pour tout $\sigma \in \mathfrak{S}_p$ (o\`u $\varepsilon(\sigma)$ d\'esigne  la signature de $\sigma$),
et
\begin{multline*}
\quad P(X_1,\ldots, X_{j-1}, -X_j-n, X_{j+1}, \ldots, X_p )
\\
= (-1)^{A(n+1)+1} P( X_1,\ldots, X_{j-1}, X_j, X_{j+1}, \ldots, X_p )\quad
\end{multline*}
pour tout $j \in \{1, \ldots, p\}.$ On suppose que $P$ est de degr\'e au plus $A(n+1)-2$
par rapport \`a chacune des variables. Alors la s\'erie
$$
\sum_{k_1\ge \cdots \ge k_p\ge 1}
\frac{P(k_1,  \ldots, k_p)}{(k_1)_{n+1}^{A} \cdots (k_p)_{n+1}^{A}}
$$
est convergente et c'est un polyn\^ome \`a coefficients rationnels en les quantit\'es
\begin{equation} \label{eqzetaas}
\sum_{\sigma\in\mathfrak{S}_q} \varepsilon(\sigma) \,\zeta(s_{\sigma(1)}, \ldots, s_{\sigma(q)})
\end{equation}
avec $q \in  \{1, \ldots, p\}$ et $s_1, \ldots, s_q \geq 3$ impairs.
\end{theo}
La somme \eqref{eqzetaas} est appel\'ee  {\em polyz\^eta antisym\'etrique} dans \cite{CFRsym}.
Lorsque $q=1$, il s'agit simplement de $\zeta(s_1)$. Pour $q=2$, on obtient $\zeta(s_1, s_2) - \zeta(s_2, s_1)$.
L'\'enonc\'e plus pr\'ecis donn\'e dans \cite{CFRsym} montre notamment que lorsque $p=1$,
on obtient une forme lin\'eaire en 1 et les $\zeta(s)$, pour $s$ impair compris entre 3
et $A$. Quand $p=2$, on obtient une forme lin\'eaire en 1, les $\zeta(s)$  pour $s$
impair compris entre 3 et $2A$, et les  $\zeta(s, s') - \zeta(s',s)$  pour $s, s'$ impairs tels que $3 \leq s < s' \leq A$.
\footnote{Les trois premiers exemples num\'eriques pr\'ec\'edant le th\'eor\`eme~\ref{conj:anti-sym}
sugg\`erent que, pour $p=2$, on a parfois des
z\^etas simples jusqu'\`a $2A-3$ et des z\^etas doubles avec $3\le s<s'\le A-2$ seulement.
Nous n'avons pas cherch\'e \`a savoir sous quelles conditions cela est vrai.}

Enfin, une derni\`ere difficult\'e, et non la moindre,
consiste \`a obtenir
des r\'esultats diophantiens en direction des
conjectures~\ref{conj:KZ} et~\ref{conj:zeta alg inde} \`a l'aide de l'approche
combinatoire d\'evelopp\'ee ici. Nous nous contentons
ici de d\'emontrer un th\'eor\`eme
{\og technique \fg} concernant le d\'enomi\-na\-teur commun aux coefficients
rationnels des combinaisons lin\'eaires
produites par certaines s\'eries du type
de~\eqref{eq:serie mult generale} : voir le th\'eor\`eme~\ref{theo:raffinement} au
paragraphe~\ref{sec:raffinement algo}.

\bigskip

\noindent{\bf Remerciements : } Nous avons pu faire fonctionner notre impl\'ementation de l'algorithme pr\'esent\'e dans ce texte sur la grappe M\'edicis. Cela nous a permis de diminuer les temps de calcul n\'ecessaires.

\section{Liens avec les int\'egrales hyperg\'eom\'etriques}

Dans ce paragraphe, on s'int\'eresse au lien entre certaines
int\'egrales multiples naturellement li\'ees aux polyz\^etas
et les s\'eries multiples que nous consid\'erons dans le pr\'esent article.
 \`A nos yeux, la souplesse combinatoire des s\'eries semble bien adapt\'ee \`a la construction de formes
lin\'eaires en polyz\^etas mais l'utilisation d'une int\'egrale ou d'une s\'erie dans
ce but est essentiellement une affaire de go\^ut, chacune
ayant des avantages et des inconv\'enients. De plus, nous mentionnons
certaines int\'egrales dont on sait qu'elles s'expriment \`a  l'aide de polyz\^etas mais auxquelles nos m\'ethodes ne s'appliquent pas.

\subsection{Exemples}

Il n'est pas possible de citer l'ensemble des int\'egrales multiples hyperg\'eom\'etriques qui sont
apparues dans la litt\'erature et nous ne mentionnons que les exemples les plus
connus.

Posons, pour tous entiers $A\ge 2$ et $n\ge 0$,
$$
J_{A,n}=\int_{[0,1]^A}\frac{\prod_{j=1}^A x_j^n(1-x_j)^n}
{Q_A(x_1,x_2,\ldots, x_A)^{n+1}}
\,\dd x_1\cdots \dd x_A,
$$
o{\`u} $Q_A(\underline x)=
1-(\cdots (1-(1-x_A)x_{A-1})\cdots)x_1$. Lorsque $A=2$ et $A=3$, on  retrouve les c\'el\`ebres int\'egrales de
Beukers~\cite{be}, qui a red\'emontr\'e le th\'eor\`eme d'Ap\'ery
en utilisant le fait que
$
J_{2,n}
\in \mathbb{Q} +  \mathbb{Q} \zeta(2)$ et $J_{3,n}
\in
\mathbb{Q} +  \mathbb{Q} \zeta(3),
$
sous une forme plus pr\'ecise.
En restant en dimension $A=2$ ou $A=3$, ces int\'egrales ont ensuite \'et\'e g\'en\'eralis\'ees dans le
but d'am\'eliorer les mesures d'irrationalit\'e respectives de $\zeta(2)$
et $\zeta(3)$~: le point d'orgue est la {\og m\'ethode du groupe \fg} de Rhin-Viola~\cite{rv1, rv2}, qui ont suivi des travaux de
Hata~\cite{hata2, hata} en particulier. La principale difficult\'e de cette approche
consiste \`a montrer directement que ces int\'egrales sont bien des formes
lin\'eaires en les valeurs de z\^eta.

En dimension sup\'erieure, Vasilyev~\cite{vasilyev} a
 formul{\'e} la conjecture
suivante, qu'il a prouv{\'e}e pour $A=4$ et $5$ :
{\em Pour tous entiers $A\ge 2$ et $n\ge 0$, il existe
des rationnels $(p_{j,A,n})_{j=0, 2, 3, \ldots, A}$ tels que
\begin{equation*}\label{eq:vasint}
J_{A,n}=p_{0,A,n}+\sum_{j \equiv A \,(\textup{mod}\, 2)} p_{j,A,n} \zeta(j).
\end{equation*}}
Cette conjecture, dont l'attaque directe est tr\`es difficile,  a {\'e}t{\'e} d{\'e}montr{\'e}e
par Zudilin \cite[paragraphe~8]{zucaen} au moyen
d'une identit{\'e} inattendue entre les int{\'e}grales
de Vasilyev et certaines  s{\'e}ries hyperg{\'e}om{\'e}triques
tr{\`e}s bien {\'e}quilibr{\'e}es. Comme on le montre au paragraphe~\ref{ssec:serieunevar}, il est alors assez facile d'obtenir une forme
lin\'eaire en
valeurs de z\^eta \`a partir d'une s\'erie hyperg\'eom\'etrique simple.

Il existe par ailleurs des int\'egrales d'une forme assez diff\'erente et qui ont \'et\'e \'etudi\'ees principalement par
Sorokin~\cite{sorokin2, so1}.
Dans~\cite{so1}, il a obtenu une preuve alternative du th\'eor\`eme d'Ap\'ery en montrant que
\begin{equation}\label{eq:intesorokinzeta3}
S_{3,n}=\int_{[0,1]^3} \frac{x^n(1-x)^ny^n(1-y)^nz^n(1-z)^n}{(1-xy)^{n+1}(1-xyz)^{n+1}}
\, \dd x\dd y \dd z\in
\mathbb{Q} +  \mathbb{Q} \zeta(3),
\end{equation}
tandis que dans~\cite{sorokin2},
il a obtenu une nouvelle preuve de la transcendance de $\pi$ en utilisant
l'int\'egrale
\begin{equation}\label{eq:padepisorokin}
T_{A,n}=\int_{[0,1]^{2A}} \prod_{j=1}^A \frac{(x_jy_j)^{n+(A-j)(n+1)}
(1-x_j)^n(1-y_j)^n}{(1-x_1y_1\cdots x_jy_j)^{n+1}} \, \dd x_j \dd y_j,
\end{equation}
dont il a montr\'e qu'elle
\'etait une forme lin\'eaire rationnelle
en 1 et les $\zeta(2,2,\ldots, 2)=\pi^{2j}/(2j+1)!$, pour $j=1, \ldots, A$,
lorsque $z=1$. Dans les deux cas,
Sorokin parvient \`a exprimer ses int\'egrales comme combinaison lin\'eaire de valeurs de polyz\^etas en r\'esolvant de mani\`ere
it\'erative des
probl\`emes de Pad\'e non triviaux. D'une mani\`ere g\'en\'erale, lorsqu'une int\'egrale provient
d'un probl\`eme de Pad\'e explicite, il arrive que l'\'enonc\'e m\^eme du probl\`eme
permettent d'\'eliminer
{\em a priori } certains polyz\^etas des formes lin\'eaires lorsque l'on sp\'ecialise les polylogarithmes multiples
en $1$ ou autre valeur int\'eressante. Ceci
conf\`ere un grand avantage \`a cette approche lorsqu'on peut la mettre en {\oe}uvre mais elle
semble difficile \`a g\'en\'eraliser. De fait,
les travaux ult\'erieurs cherchent tous \`a s'affranchir de l'\'etape {\og Pad\'e\fg}.

Le fait particuli\`erement remarquable que $J_{3,n}=S_{3,n}$ pour tout entier $n\ge 0$ a \'et\'e g\'en\'eralis\'e par
Fischler~\cite{fi2}
et Zlobin~\cite{zlo} ind\'ependamment, qui ont montr\'e entre autres choses que l'on a
les identit\'es
\begin{equation}\label{eq:vasi=soro}
J_{A,n}= \int_{[0,1]^A}\prod_{j=1}^{A/2} \frac{x_j^n (1-x_j)^ny_j^n(1-y_j)^n}{(1-x_1y_1\cdots x_jy_j)^{n+1}}
\, \dd x_j \dd y_j
\end{equation}
pour $A\ge 2$ pair et
\begin{equation}\label{eq:vasi=soro2}
J_{A,n}= \int_{[0,1]^A}\frac{z^n(1-z)^n}{(1-x_1y_1\cdots x_{a} y_{a} z)^{n+1}}
\bigg(\prod_{j=1}^{a} \frac{x_j^n (1-x_j)^ny_j^n(1-y_j)^n}{(1-x_1y_1\cdots x_jy_j)^{n+1}}
\, \dd x_j \dd y_j \bigg)\dd z
\end{equation}
pour $A\ge 3$ impair avec $a=(A-1)/2$.
Il d\'ecoule de ces travaux l'intuition assez nette que l'on ne perd rien \`a travailler avec des g\'en\'eralisations
de l'une ou l'autre des int\'egrales
$J_{A,n}$ et $S_{A,n}$. Il
s'av\`ere que les int\'egrales de Sorokin $S_{A,n}$ \`a droite de~\eqref{eq:vasi=soro} et~\eqref{eq:vasi=soro2}
se d\'eveloppent un peu plus facilement en s\'eries multiples que les int\'egrales $J_{A,n}$ et qu'elles
donnent imm\'ediatement  des polyz\^etas dans le cas $n=0$. Dans une perspective diophantienne, il est donc
naturel de produire des formes lin\'eaires en polyz\^etas \`a partir d'int\'egrales du type de Sorokin les plus g\'en\'erales
possibles ; une telle relation  a \'et\'e d\'emontr\'ee par Zlobin~\cite{zlo}. La proposition~\ref{prop:intesoro=seriemult} (d\'emontr\'ee
au paragraphe~\ref{sec:prop:intesoro=seriemult} ci-dessous) coupl\'ee
aux r\'esultats de cet article nous permet de red\'emontrer une assertion similaire \`a celle de Zlobin mais nous insistons ici sur le
fait que nos r\'esultats (r\'esum\'es informellement par le th\'eor\`eme~\ref{theo:informel})
nous permettent de traiter des s\'eries multiples plus g\'en\'erales que celles apparaissant dans la
proposition~\ref{prop:intesoro=seriemult} ou dans les travaux de
Zlobin.

Terminons ce paragraphe en mentionnant un r\'ecent article de Zlobin~\cite{zlobin2}, o\`u il obtient
une nouvelle preuve de la conjecture de Vasilyev en partant de
l'int\'egrale $S_{A,n}$ convenablement d\'evelopp\'ee en s\'erie multiple~: il
s'agit d'un remarquable tour de force.

\subsection{D\'eveloppement en s\'erie de certaines int\'egrales de Sorokin} \label{sec:prop:intesoro=seriemult}

Le but de ce paragraphe est d'exprimer une int\'egrale de type Sorokin relativement g\'en\'erale
(elle contient du moins tous les cas mentionn\'es ci-dessus) comme une s\'erie multiple. Cette derni\`ere est un cas particulier de celle
que nous d\'eveloppons en
polylogarithmes multiples et/ou polyz\^etas dans la suite de l'article : pour ceux
qui aiment travailler \`a partir d'int\'egrales,
la proposition~\ref{prop:intesoro=seriemult} (voir aussi le lemme 2 de \cite{ZlobinZametki2005}) est donc la premi\`ere \'etape de notre algorithme de construction de formes lin\'eaires en
polyz\^etas.

\begin{prop} \label{prop:intesoro=seriemult}
Soient des entiers $D, p\ge 1$ et des entiers positifs
$r_1, \ldots, r_p$, $s_1, \ldots, s_p$,  $t_1, \ldots, t_p$ et $0=d_0< d_1 < d_2<\cdots <d_p=D$. Pour
tout  complexe $z$ tel que $ \vert z\vert >1$, on~a l'identit\'e
\begin{multline}\label{eq:prop:intesoro=seriemult}
\int_{[0,1]^D} \prod_{j=1}^p \frac{
\prod_{\ell=d_{j-1}+1}^{d_j}x_\ell^{r_j}(1-x_\ell)^{s_j}}
{(z-x_1\cdots x_{d_j})^{t_j+1}}\, \dd x_j
\\
= z^{-(t_1+\cdots+ t_p + p-1)} \cdot \prod_{j=1}^{p}\frac{s_j!^{A_j}}{t_j!} \cdot
\sum_{k_1 \ge \cdots \ge k_p\ge 1}  z^{-k_1} \prod_{j=1}^p\frac{(k_j-k_{j+1}+1)_{t_j}}
{(k_j+r_j)_{s_j+1}^{A_j}},
\end{multline}
o\`u $k_{p+1}=1$ et $A_j=d_j-d_{j-1}$ pour $j=1, \ldots, p.$ La s\'erie
est de profondeur $p$ et de poids~$D$.
\end{prop}

\begin{remark} L'\'equation~\eqref{eq:prop:intesoro=seriemult} s'\'etend
\`a $\vert z\vert =1$ lorsque les deux membres ont un sens simultan\'ement.

Dans les applications diophantiennes, il est pratique de sommer sur des indices
$K_j$ d\'efinis par $K_j=k_j+r$, o\`u $r= \min r_j$. En particulier, si tous les $r_j$ sont
\'egaux \`a $r$, la s\'erie s'\'ecrit
$$
\sum_{K_1 \ge \cdots \ge K_p\ge r+1} z^{-K_1-r} \prod_{j=1}^p\frac{(K_j-K_{j+1}+1)_{t_j}}
{(K_j)_{s_j+1}^{A_j}}
$$
avec $K_{p+1}=r+1$. De plus, si $t_p=r$,
la pr\'esence du symbole de Pochhammer
$(K_p-r)_{r}$ implique que sommer sur l'ensemble d'indices
$K_1 \ge \cdots \ge K_p\ge r+1$  revient au m\^eme que sommer sur $K_1 \ge \cdots \ge K_p\ge 1$.
\end{remark}

\begin{proof} Supposer que $\vert z\vert >1$ assure que les diverses d\'ecompositions en s\'eries et
 inversions s\'eries-int\'egrales ci-dessous sont licites. Le cas
 d'un point du cercle $\vert z\vert =1$
s'obtient  en invoquant des crit\`eres de continuit\'e (th\'eor\`emes d'Abel, de Lebesgue, etc).

On d\'eveloppe le d\'enominateur de l'int\'egrale multiple, not\'ee $I(z)$ dans la suite, au moyen de
l'identit\'e (avec $\vert z\vert >1, 0\le x\le 1$) :
$$
\frac{1}{(z-x)^{t+1}} = \frac{1}{z^t} \sum_{m=0}^{\infty} \binom{m+t}{m}
\bigg(\frac{x}{z}\bigg)^m
$$
et on obtient alors
\begin{multline*}
I(z)
= z^{-(t_1+\cdots+ t_p + p)}
\\\times
\sum_{m_1, \ldots, m_p \ge 0} \prod_{j=1}^p \binom{m_j+t_j}{m_j} z^{-m_j} \int \limits_{[0,1]^D}
\prod_{j=1}^p \bigg((x_1\cdots x_{d_j})^{m_j} \prod_{\ell=d_{j-1}+1}^{d_j}x_\ell^{r_j}(1-x_\ell)^{s_j} \,\dd x_j\bigg).
\end{multline*}
Or on v\'erifie que
$$
\prod_{j=1}^p \bigg((x_1\cdots x_{d_j})^{m_j} \prod_{\ell=d_{j-1}+1}^{d_j}x_\ell^{r_j}(1-x_\ell)^{s_j}
\bigg) = \prod_{j=1}^p   \bigg(\prod_{\ell=d_{j-1}+1}^{d_j} x_\ell^{r_j+m_j+\cdots + m_p}
(1-x_\ell)^{s_j}\bigg).
$$
On peut s\'eparer les variables dans l'int\'egrale et on obtient alors $D$
int\'egrales facilement calculables (ce sont des
fonctions Beta d'Euler), d'o\`u
$$
I(z)= z^{-(t_1+\cdots +t_p + p)}
\sum_{m_1, \ldots, m_p \ge 0} \prod_{j=1}^p
\frac{z^{-m_j} \binom{m_j+t_j}{m_j}}{\binom{r_j+s_j+m_j+\cdots +m_p}{s_j}^{A_j}(r_j+s_j+m_j+\cdots +m_p+1)^{A_j}}.
$$
On utilise maintenant les deux transformations triviales
$
\binom{m_j+t_j}{m_j}  = \frac{(m_j+1)_{t_j}}{t_j!}
$
et
\begin{multline*}
\qquad \binom{r_j+s_j+m_j+\cdots +m_p}{s_j}(r_j+s_j+m_j+\cdots +m_p+1)
\\
 = \frac{(r_j+m_j+\cdots +m_p+1)_{s_j+1}}{s_j!}\qquad
\end{multline*}
et on pose $k_j=m_j+\cdots +m_p+1$ pour $j=1, \ldots, p$, ainsi que $k_{p+1}=1$. On obtient alors
$$
I(z)= z^{-(t_1+\cdots+ t_p + p-1)} \cdot \prod_{j=1}^{p}\frac{s_j!^{A_j}}{t_j!}\cdot
\sum_{k_1 \ge  \cdots \ge k_p \ge 1}  z^{-k_1} \prod_{j=1}^p\frac{(k_j-k_{j+1}+1)_{t_j}}
{(k_j+r_j)_{s_j+1}^{A_j}},
$$
ce qui termine la preuve.
\end{proof}

\`A titre d'exemples, remarquons que
l'int\'egrale $S_{3,n}$ en~\eqref{eq:intesorokinzeta3}
vaut exactement la s\'erie~\eqref{eq:sorokin1} donn\'ee dans l'introduction tandis que
l'int\'egrale~\eqref{eq:padepisorokin} s'exprime de la mani\`ere suivante, en posant
$k_{A+1} = n+1$ :
$$
T_{A,n}= n!^{A} \sum_{k_1 \ge \cdots \ge k_A\ge 1} \prod_{j=1}^A
\frac{(k_j -k_{j+1}+1)_n}{(k_j+(A-j)(n+1))_{n+1}^2}.
$$
Notre algorithme permet ensuite d'exprimer effectivement ces int\'egrales comme
des formes lin\'eaires en polyz\^etas.  Jusqu'en poids 4, on ne voit appara\^itre que des valeurs de z\^eta
car tous les polyz\^etas de poids $\le 4$ sont des multiples
rationnels de 1, $\zeta(2)$, $\zeta(3)$ ou $\zeta(4)$. En revanche, \`a partir du poids 5, on doit s'attendre \`a obtenir des
polyz\^etas lin\'eairement ind\'ependants (du moins, conjecturalement) des valeurs de z\^eta, comme le montre l'exemple
de l'int\'egrale
$$
\int_{[0,1]^5} \frac{ \prod_{j=1}^5 x_j^n(1-x_j)^n dx_j}
{ (1-x_1x_2x_3)^{n+1}(1-x_1x_2x_3x_4x_5)^{n+1} }.
$$
Pour $n=0$, elle vaut $\zeta(3,2) = -11 \zeta(5) /2 +3 \zeta(2)\zeta(3)$, qui n'est
donc probablement pas un multiple rationnel d'une valeur de z\^eta en un
 entier. Pour $n=1$, $2$ et $3$, elle est une combinaison lin\'eaire rationnelle en
$\zeta(2), \zeta(3), \zeta(4), \zeta(5), \zeta(2,2)$ et $\zeta(3,2)$. Le coefficient de
$\zeta(3,2)$ dans ces combinaisons lin\'eaires est non nul ; on peut l'expliciter comme une somme double finie, ce qui pourrait peut-\^etre permettre de d\'emontrer
  qu'il est non nul pour tout $n \ge 4$, si n\'ecessaire.

\subsection{D'autres exemples d'int\'egrales hyperg\'eom\'etriques}

Il existe beaucoup d'autres types d'int\'egrales hyperg\'eom\'etriques
que celles de Vasilyev et Sorokin et dont par des moyens
plus ou moins d\'etourn\'es on sait qu'elles s'expriment comme
formes lin\'eaires en polyz\^etas. Les deux exemples que nous allons aborder sont
dus \`a Zudilin et Goncharov-Manin respectivement.

L'int\'egrale consid\'er\'ee par Zudilin est la suivante :
$$
Z_n= \int_{[0,1]^5}
\frac{\prod_{j=1}^5 x_j^n(1-x_j)^n \,\dd x_j}{Q(x_1,x_2,x_3,x_4,x_5)^{n+1}}
$$
o\`u
$
Q(\underline x)=x_1(1-(1-(1-(1-x_2)x_3)x_4)x_5)+ (1-x_1x_2x_3x_4x_5).
$
Par un proc\'ed\'e indirect (bas\'e sur des transformations hyperg\'eom\'etriques),
il montre que $Z_n$ est \'egale \`a une
s\'erie de nature hyperg\'eom\'etrique tr\`es bien \'equilibr\'ee (avec double une
d\'erivation du sommande) et il en d\'eduit que
$Z_n\in\mathbb{Q}+\mathbb{Q}\zeta(4)$.

Les int\'egrales de Goncharov-Manin~\cite{goncharov2}
apparaissent quant \`a elles comme des
p\'eriodes de certains motifs de Tate mixte, dont Brown~\cite{brown} a donn\'e la forme explicite
suivante:
\begin{equation}\label{eq:gonchmanin}
\int_{[0,1]^A} \frac{\prod_{j=1}^A x_j^{r_j}(1-x_j)^{s_j}\,\dd x_j}
{\prod_{1\le i<j\le A} (1-x_i\cdots x_j)
^{t_{i,j}}}
\end{equation}
avec des entiers $r_j, s_j, t_{i,j}\ge 0$ tels que l'int\'egrale converge.
Remarquons que~\eqref{eq:gonchmanin} contient comme cas particulier les int\'egrales
abord\'ees par la proposition~\ref{prop:intesoro=seriemult} (en $z=1$).
Par des arguments de nature g\'eom\'etrique, Brown a prouv\'e une conjecture de
Goncharov-Manin qui affirmait que ces int\'egrales sont toujours de formes lin\'eaires rationnelles en polyz\^etas. Sa m\'ethode n'est
malheureusement pas constructive, ce qui rend impossible une quelconque utilisation diophantienne de son th\'eor\`eme par les voies
classiques.

Ces deux types d'int\'egrales ont donc le d\'efaut de n'\^etre \'evaluable que par des
proc\'ed\'es tr\`es indirects. Pour rem\'edier \`a cela, on pourrait tenter de les
{\og d\'evelopper \fg}  en s\'eries multiples \`a la mani\`ere de la
proposition~\ref{prop:intesoro=seriemult}, puis esp\'erer appliquer une g\'en\'eralisation convenable de notre algorithme.
Ceci n'aura rien
d'\'evident ;
par exemple, les cas les plus simples de l'int\'egrale~\eqref{eq:gonchmanin}
peuvent conduire \`a des s\'eries telles que
$$
\sum_{m, n\ge 1} \frac{1}{m^{s_1}n^{s_2}(m+n)^{s_3}},
$$
dont il n'est m\^eme pas clair qu'elles
puissent s'exprimer \`a l'aide de polyz\^etas (c'est cependant
bien le cas : voir~\cite{moll} pour plus de d\'etails et des r\'ef\'erences).  \'Etendre notre algorithme n\'ecessitera donc des id\'ees
nouvelles.

\section{\'Etude de deux situations instructives}\label{sec:prof 1 et 2}

\subsection{Le cas de la profondeur 1}\label{ssec:serieunevar}
La strat\'egie~\footnote{C'est essentiellement la seule dont on dispose~: toutes les autres
approches connues produisent les m\^emes formes lin\'eaires (voir~\cite{fi}).}
pour d\'emontrer les  th\'eor\`emes diophantiens concernant les valeurs de la  fonction z\^eta est la suivante.
Soient des entiers $n\ge 0$, $A\ge 1$ et $P(X)\in\mathbb{Q}[X]$. Consid\'erons
la fraction rationnelle
$\dis R(X)= P(X)/(X)_{n+1}^A$ ainsi que
la s\'erie
$$
S(z) = \sum_{k=1}^{\infty} R(k) \,z^{-k}.
$$
On suppose cette derni\`ere convergente pour $z= 1$, ce qui impose que
deg$(P)$ $\le A(n+1)-2$. On commence par d\'evelopper $R(X)$
en \'el\'ements simples :
$$
R(X) = \sum_{s=1} ^A\sum_{j=0}^n
\frac{C\bigg[\,\begin{matrix} s\\j\end{matrix}\,\bigg]}{(X+j) ^s} \quad
\textup{avec} \quad
C\bigg[\,\begin{matrix} s\\j\end{matrix}\,\bigg]=
\frac{1}{(A-s)!}\bigg(R(X)(X+j)^A\bigg) ^{(A-s)}\bigg\vert_{X=-j}
$$
et, en reportant dans $S(z)$, on obtient
\begin{eqnarray*}
S(z)= \sum_{s=1} ^A\sum_{j=0}^n
C\bigg[\,\begin{matrix} s\\j\end{matrix}\,\bigg] \sum_{k=1}^{\infty}
\frac{z ^{-k}}{(k+j)^s}.
\end{eqnarray*}

On remarque alors que, trivialement,
\begin{equation}\label{eq:brique prof 1}
\sum_{k=1}^{\infty} \frac{z ^{-k}}{(k+j)^s} = z^j \,\li_s(1/z) -
\sum_{k=1}^{j} \frac{z ^{j-k}}{k^s}
\end{equation}
et donc qu'il existe des polyn\^omes
$Q_{s}(z)$ $\in\mathbb{Q}[z]$, de degr\'e au plus $n$,
tels que
$$
S(z)=Q_{0}(z)+ \sum_{s=1}^A Q_{s}(z)\Li_{s}(1/z).
$$
On a bien s\^ur
$\li_s(1)=\zeta(s)$ et
$\li_s(-1)= (2^{1-s}-1)\zeta(s)$ pour tout $s>1$. Pour $s\ge 1$, on a l'expression tr\`es simple
$$
Q_s(z) =\sum_{j=0}^n C\bigg[\,\begin{matrix} s\\j\end{matrix}\,\bigg] z^j.
$$

Pour les applications envisag\'ees, il est important de se ramener
\`a des coefficients entiers et on montre
que $Q_{1}(1)=0$ et
 $\dd_n^{A-j} \,Q_{j}(z)\in\mathbb{Z}[z]$ pour tout $j\in\{0, \ldots, A\}$,
o\`u  $\dd_n=\textup{p.p.c.m.}\{1, 2, \ldots, n\}$.
Il existe donc des entiers $q_{j}$ tels que
$$
\dd_n^{A}S(1) = q_{0} + \sum_{s=2}^{A} q_{s} \,\zeta(s),
$$
et une expression similaire pour  $S(-1)$.

Tout le probl\`eme r\'eside  maintenant dans des choix
de $A$ et de $P$ tels que l'on puisse appliquer efficacement
un crit\`ere d'irrationalit\'e ou d'ind\'e\-pen\-dan\-ce lin\'eaire :
il appara{\^\i}t rapidement que l'on
doit \'eliminer les nombres $\zeta(s)$ pour $s$ pair, sous peine de n'obtenir que des
r\'esultats triviaux. Une mani\`ere d'y parvenir est d'imposer que le polyn\^ome
 $P(X)$ satisfasse \`a
\begin{equation}\label{eq:symetrie prof 1}
P(-X-n)=-P(X).
\end{equation}
En effet, par unicit\'e de la d\'ecomposition de
$R(X)$ en \'el\'ements simples, l'\'equation~\eqref{eq:symetrie prof 1}
se traduit par
$C\bigg[\,\begin{matrix} s\\n-j\end{matrix}\,\bigg]
=(-1)^{A(n+1)+s+1}
C\bigg[\,\begin{matrix} s\\j\end{matrix}\,\bigg]$ et donc
les coefficients $q_{s}$ sont nuls pour $s$ pair lorsque $A$ est lui-m\^eme pair.
Par exemple, lorsque $A$ est pair, on peut utiliser les s\'eries
\begin{equation*}
n!^{A-2r}
\sum_{k=1}^{\infty}
\left(k+\frac{n}{2}\right)\frac{(k-rn)_{rn}(k+n+1)_{rn}}{(k)_{n+1}^A}
=q_{0}+\sum_{{s=3 \atop s \,\textup{impair}}}^A q_{s} \,\zeta(s),
\end{equation*}
qui sont des s\'eries hyperg\'eom\'etriques sp\'eciales, dites
{\em very-well-poised} (voir~\cite{ba, sl} pour la d\'efinition exacte).
On se r\'ef\`erera \`a~\cite{br, fi, kratriv, ri, zud} pour
plus de d\'etails sur l'utilisation diophantienne de ce type de s\'erie.

\subsection{Le cas de la profondeur 2} \label{ssec: algo prof 2}

Une fois formalis\'e le cas de la profondeur 1, il est naturel d'essayer de suivre la m\^eme
d\'emarche en profondeur sup\'erieure. Le cas de la profondeur $p=2$ est d\'ej\`a instructif
et nous allons le traiter en d\'etails.

Nous expliquons notre approche sur la s\'erie
suivante
\begin{equation}
\label{eq:sorokin2}
S(z_1,z_2)=\sum_{k_1\ge k_2\ge 1}
\frac{P(k_1,k_2)}{(k_1)_{n+1}^2(k_2)_{n+1}^2}
\,z_1^{-k_1}z_2^{-k_2},
\end{equation}
avec $\textup{deg}_{k_1}(P)\le A(n+1)-2$ et
$\textup{deg}_{k_2}(P)\le A(n+1)-2$, ce qui assure que la s\'erie
converge absolument pour $\vert z_1\vert \ge 1$ et
$\vert z_2\vert \ge 1$. On notera que
la s\'erie introduite en~\eqref{eq:sorokin1} ne v\'erifie pas cette
condition de degr\'e : les cons\'equences de cela sont \'evoqu\'ees \`a
la fin de ce paragraphe.

La premi\`ere \'etape consiste, comme pr\'ec\'edemment, \`a d\'ecompo\-ser
en \'el\'ements simples la fraction rationnelle
qui constitue le sommande de
$S(z_1,z_2)$  :
$$
\frac{P(k_1,k_2)}{(k_1)_{n+1}^2(k_2)_{n+1}^2}
= \sum_{j_1, j_2=0}^n
  \sum_{s_1, s_2=1}^2
\frac{C\bigg[\,\begin{matrix} s_1, s_2\\j_1, j_2\end{matrix}\,\bigg]}{(k_1+j_1)^{s_1}(k_2+j_2)^{s_2}},
$$
o\`u les $C\bigg[\,\begin{matrix} s_1, s_2\\j_1, j_2\end{matrix}\,\bigg]$
sont des rationnels explicitables.
Il est important de noter que la condition portant sur les
degr\'es de $P$ implique que cette d\'ecomposition n'a pas de partie enti\`ere.
En reportant dans $S(z_1,z_2)$, on obtient ainsi
\begin{equation*}
S(z_1,z_2)=
\sum_{j_1, j_2=0}^n \, \sum_{s_1, s_2=1}^2
C\bigg[\,\begin{matrix} s_1, s_2\\j_1, j_2\end{matrix}\,\bigg]
\sum_{k_1\ge k_2\ge 1} \frac{z_1^{-k_1}z_2^{-k_2}}
{(k_1+j_1)^{s_1}(k_2+j_2)^{s_2}}.
\end{equation*}

La deuxi\`eme \'etape consiste \`a exprimer
explicitement la s\'erie
\begin{equation}\label{eq:brique double elementaire}
\sum_{k_1 =1}^{\infty}
\frac{z_1^{-k_1}}{(k_1+j_1)^{s_1}}\,\sum_{k_2=1}^{k_1} \frac{z_2^{-k_2}}{(k_2+j_2)^{s_2}}
\end{equation}
comme une combinaison lin\'eaire \`a
coefficients dans
$\mathbb{Q}[z_1^{\pm 1}, z_2^{\pm 1}]$ en
les polylogarithmes multiples (larges ou stricts).
Comme on l'a vu en~\eqref{eq:brique prof 1}, dans le cas d'une seule variable ($p=1$), c'est une
\'etape triviale mais, malheureusement, en deux variables, ce n'est plus le cas.
On \'ecrit tout d'abord la somme int\'erieure sur $k_2$ comme
\begin{equation*}
\sum_{k_2=1}^{k_1}\frac{z_2^{-k_2}}{(k_2+j_2)^{s_2}}
=\sum_{k_2=j_2+1}^{k_1+j_2}\frac{z_2^{j_2-k_2}}{k_2^{s_2}}
=
\bigg(\sum_{k_2=1}^{k_1+j_1}-\sum_{k_2=1}^{j_2}+\;\varepsilon_{j_1,j_2}
\sum_{k_2=k_1+j_1 \wedge j_2+1}^{k_1+j_1 \vee j_2}\bigg)\frac{z_2^{j_2-k_2}}{k_2^{s_2}}
\end{equation*}
o\`u $j_1 \wedge j_2 = \min(j_1,j_2)$,
$j_1\vee j_2 =\max(j_1, j_2)$ et
$\varepsilon_{j_1,j_2}=1$ si $j_1<j_2$, $-1$ si
$j_1>j_2$, $0$ si $j_1=j_2$.  Puis on
reporte ces trois sommes dans la somme sur $k_1$.
Les deux premi\`eres s\'eries se traitent facilement~:
\begin{eqnarray*}
\sum_{k_1=1}^{\infty}\frac{z_1^{-k_1}}{(k_1+j_1)^{s_1}}
\sum_{k_2=1}^{k_1+j_1}\frac{z_2^{j_2-k_2}}{k_2^{s_2}}
 & = & \sum_{k_1=j_1+1}^{\infty}\frac{z_1^{j_1-k_1}}{k_1^{s_1}}
 \sum_{k_2=1}^{k_1}\frac{z_2^{j_2-k_2}}{k_2^{s_2}}\\
 & = & z_1^{j_1}z_2^{j_2}\;\textup{La}_{{s_1},{s_2}}(1/z_1,1/z_2)-
\sum_{k_1=1}^{j_1}\frac{z_1^{j_1-k_1}}{k_1^{s_1}}
\sum_{k_2=1}^{k_1}\frac{z_2^{j_2-k_2}}{k_2^{s_2}}
\end{eqnarray*}
et
\begin{multline*}
\sum_{k_1=1}^{\infty}\frac{z_1^{-k_1}}{(k_1+j_1)^{s_1}}
\sum_{k_2=1}^{j_2}\frac{z_2^{j_2-k_2}}{k_2^{s_2}}
\\
=z_1^{j_1}\bigg(\sum_{k_2=1}^{j_2}
\frac{z_2^{j_2-k_2}}{k_2^{s_2}}\bigg)\textup{La}_{{s_1}}(1/z_1)-
\bigg(\sum_{k_1=1}^{j_1}\frac{z_1^{j_1-k_1}}{k_1^{s_1}}\bigg)
\bigg(\sum_{k_2=1}^{j_2}\frac{z_2^{j_2-k_2}}{k_2^{s_2}}\bigg).\qquad
\end{multline*}
La troisi\`eme s\'erie est un peu plus compliqu\'ee : on note que
\begin{eqnarray*}
\sum_{k_1=1}^{\infty}\frac{z_1^{-k_1}}{(k_1+j_1)^{s_1}}
\sum_{k_2=k_1+j_1 \wedge j_2+1}^{k_1+j_1 \vee j_2}\frac{z_2^{j_2-k_2}}{k_2^{s_2}}
 & = &
\sum_{k_2=j_1 \wedge j_2+1}^{j_1 \vee j_2} z_2^{j_2-k_2}
\sum_{k_1=1}^{\infty}
\frac{(z_1z_2)^{-k_1}}{(k_1+j_1)^{s_1}(k_1+k_2)^{s_2}}
\end{eqnarray*}
puis l'on d\'eveloppe en \'el\'ements simples la fraction rationnelle
$$
\frac{1}{(k_1+j_1)^{s_1}(k_1+k_2)^{s_2}}
$$
pour conclure
que cette s\'erie s'\'ecrit comme une combinaison lin\'eaire
de $\textup{La}_s(1/z_1z_2)$ avec $1\le s\le s_1 \vee s_2$ et aussi
de $\textup{La}_{{s_1}+{s_2}}(1/z_1z_2)$ si
$k_2=j_1 \vee j_2=j_1$,  avec des coefficients polynomiaux en $z_1^{\pm 1}$ et $z_2^{\pm 1}$.
En r\'esum\'e,  lorsque  $z_1=z_2=1$, la d\'ecomposition de la s\'erie~\eqref{eq:brique double elementaire}
fait appara{\^\i}tre au plus les polyz\^etas suivants :
$\zeta(s_1, s_2)$, $\zeta(s_1 + s_2)$ et les $\zeta(s)$ pour
$1 \le s \le s_1 \vee s_2$. En particulier, il n'y a aucune raison apparente pour que les valeurs de z\^eta
aux entiers pairs n'apparaissent pas.

Enfin, troisi\`eme \'etape, en reportant la d\'ecomposition ainsi obtenue
dans~\eqref{eq:sorokin2},
on doit identifier les polyz\^etas
qui apparaissent r\'eellement dans $S(1,1)$,
c'est-\`a-dire ceux affect\'es d'un
coefficient non-nul. Or cette identification
 n'est pas \'evidente : la s\'erie~$S(1,1)$ fait
 appara{\^\i}tre {\em a priori} les polyz\^etas
$$
\zeta(1), \,\zeta(1,1), \,\zeta(2), \,\zeta(2,1), \,\zeta(1,2),
\,\zeta(3), \,\zeta(2,2), \,\zeta(4)
$$
(certains sont
divergents). Lorsque, par exemple, $P(X_1, X_2)=(X_2-n)_n(X_1-X_2+1)_n$, il
est assez difficile de prouver que seuls $\zeta(2)$,
$\zeta(2,2)$ et $\zeta(4)$ n'ont pas un coefficient nul.

 On doit aussi parfois tenir compte
d'un autre
ph\'enom\`ene~:
contrairement \`a $S(z_1,z_2)$, la d\'ecomposition
en \'el\'ements simples du sommande de la s\'erie en~\eqref{eq:sorokin1}
produit une partie enti\`ere (puisque le degr\'e en $X_2$ de la
fraction $(X_2-n)_n(X_1-X_2+1)_n/(X_1)_{n+1}^2(X_2)_{n+1}$ est positif)
qui complique encore cette \'etape en faisant
appara{\^\i}tre des polylogarithmes multiples
{\og exotiques \fg} tels que $\la_{2,-1}(z_1,z_2)$ qu'il faut traiter de fa\c{c}on {\em ad hoc}.
Ce proc\'ed\'e   devient quasiment
inextricable en trois variables, ce qui explique
le formalisme que
nous d\'eveloppons au paragraphe~\ref{sec:l'algorithme}.


\section{D\'emonstration du th\'eor\`eme~\ref{theo:informel}}\label{sec:l'algorithme}

Nous venons de  d\'emontrer le th\'eor\`eme~\ref{theo:informel} pour $p=1$
(paragraphe~\ref{ssec:serieunevar})
et $p=2$ (paragraphe~\ref{ssec:
algo prof 2}).
Dans ce paragraphe, on le d\'emontre en toute g\'en\'eralit\'e :
la strat\'egie consiste \`a se ramener dans un premier temps \`a un cas
plus simple (paragraphe~\ref{ssec:serie hyp mult en brique}) que
l'on d\'emontre ensuite (th\'eor\`eme~\ref{prop:2} au
paragraphe~\ref{sec:2}).
Nous
en obtiendrons des raffinements au paragraphe~\ref{sec:raffinement algo}.

\subsection{D\'ecomposition des s\'eries
multiples en briques}\label{ssec:serie hyp mult en brique} \label{sec:decompositionenbrique}

En imitant le cas de la profondeur 1, nous allons transformer la
s\'erie
\begin{equation}\label{eq:S_p mult generale}
S_P\bigg[\,\begin{matrix}
A_1, \ldots, A_p\\
n_1, \ldots, n_p
\end{matrix}
\,\bigg\vert \,z_1, \ldots, z_p\, \bigg] = \sum_{k_1\ge \cdots \ge
k_p\ge 1} \frac{P(k_1,  \ldots, k_p)}{(k_1)_{n_1+1}^{A_1}
 \cdots (k_p)_{n_p+1}^{A_p}}\,z_1^{-k_1}\cdots z_p^{-k_p}
\end{equation}
en d\'eveloppant en \'el\'ements simples la fraction rationnelle
$$
R(X_1, \ldots, X_p)= \frac{P(X_1,  \ldots,
X_p)}{(X_1)_{n_1+1}^{A_1}
 \cdots (X_p)_{n_p+1}^{A_p}}.
$$
Posons $\aich=\deg_{X_i}(P)- A_i(n_i+1)$ : c'est le degr\'e en
$X_i$ de la fraction rationnelle $R$. Notons $J$ l'ensemble des
indices $i \in \unp$ tels que $\aich \geq 0$ (c'est-\`a-dire
$\deg_{X_i}(P) \geq A_i(n_i+1)$) : c'est l'ensemble des $i$ pour
lesquels $R$ est de degr\'e positif ou nul en $X_i$,
c'est-\`a-dire relativement auxquels une partie enti\`ere va
appara\^{\i}tre. Pour $I \subset \unp$, on note $\Ic = \unp \moins
I$. Alors on  a
\begin{equation} \label{eq: decomposition de R}
R(X_1, \ldots, X_p)= \sum_{I \subset J} \sum_\indtr{(s_i)_{i \in
\Ic} \mbox{ tel que }}{1 \leq s_i \leq A_i }{\mbox{pour tout } i
\in \Ic} \sum_\indtr{(j_i)_{i \in \Ic} \mbox{ tel que }}{0 \leq
j_i \leq n_i }{\mbox{pour tout } i \in \Ic} \sum_\indtr{(\sich)_{i
\in I} \mbox{ tel que }}{0 \leq \sich \leq \aich }{\mbox{pour tout
} i \in I} \Ciandco \frac{\prod_{i \in I} X_i^{\sich}}{\prod_{i
\in \Ic}(X_i+j_i)^{s_i}}
\end{equation}
avec
$$\Ciandco = \diandco \Big(R_I(Y_1, \ldots, Y_p)  \prod_{i \in \Ic} (Y_i+j_i)^{A_i}
\prod_{i \in I} Y_i^{\aich} \Big)_{\left| {\tiny \begin{array}{l} {Y_i = 0
\mbox{ pour } i \in I} \\ {Y_i = -j_i \mbox{ pour } i  \in \Ic} \end{array}}\right.}$$
en notant $\diandco$ l'op\'erateur diff\'erentiel suivant
$$
\diandco = \prod_{i \in \Ic} \Big( \frac{1}{(A_i - s_i)!} \Big(
\frac{\partial}{\partial Y_i} \Big)^{A_i-s_i} \Big)
 \prod_{i \in I}
\Big( \frac{1}{(\aich - \sich)!} \Big( \frac{\partial}{\partial
Y_i} \Big)^{\aich - \sich} \Big),
$$
et $R_I(Y_1, \ldots, Y_p)$ la fraction rationnelle obtenue \`a
partir de $R(X_1, \ldots, X_p)$ en posant :
$$
\left\{
\begin{array}{l}
X_i = \frac{\displaystyle 1}{\displaystyle Y_i} \mbox{ pour } i \in I \\
X_i = Y_i \mbox{ pour } i \in \Ic .
\end{array}
\right.
$$

Le cas particulier o\`u il n'y a pas de partie enti\`ere
correspond \`a $\aich \leq -1$ pour tout $i \in \unp$,
c'est-\`a-dire $J = \emptyset$. La somme sur $I$ se r\'eduit alors
\`a $I  = \emptyset$, la famille $(\sich)$ est vide et on obtient
la d\'ecomposition en \'el\'ements simples habituelle :
$$
R(X_1, \ldots, X_p) = \sum_{s_1=1}^{A_1}\cdots \sum_{s_p=1}^{A_p}
\sum_{j_1=0}^{n_1}\cdots \sum_{j_p=0}^{n_p}
C\bigg[\,\begin{matrix} s_1,  \ldots, s_p\\
                        j_1,  \ldots, j_p
\end{matrix} \,\bigg] \frac{1}{(X_1+j_1)^{s_1}\cdots (X_p+j_p)^{s_p} }.
$$

Revenons au cas g\'en\'eral. En reportant~\eqref{eq: decomposition
de R} dans~\eqref{eq:S_p mult generale}, on obtient
\begin{multline*}
S_P\bigg[\,\begin{matrix}
A_1, \ldots, A_p\\
n_1, \ldots, n_p
\end{matrix}
\,\bigg\vert \,z_1, \ldots, z_p\, \bigg] = \sum_{I \subset J}
\sum_\indtr{(s_i)_{i \in \Ic} \mbox{ tel que }}{1 \leq s_i \leq
A_i }{\mbox{pour tout } i \in \Ic} \sum_\indtr{(j_i)_{i \in \Ic}
\mbox{ tel que }}{0 \leq j_i \leq n_i }{\mbox{pour tout } i \in
\Ic}
\sum_\indtr{(\sich)_{i \in I}
\mbox{ tel que }}{0 \leq \sich \leq \aich }{\mbox{pour tout } i \in I}
\\
\cdot \Ciandco \sum_{k_1\ge \cdots \ge k_p\ge 1} \frac{\prod_{i \in I}
k_i^{\sich}}{\prod_{i \in \Ic}(k_i+j_i)^{s_i}} \,z_1^{-k_1}\cdots
z_p^{-k_p}.
\end{multline*}
On a donc ramen\'e le probl\`eme initial (i.e., l'\'eva\-lua\-tion de
\eqref{eq:S_p mult generale}) \`a ce\-lui de la d\'ecomposi\-tion en
polylogarithmes multiples de s\'eries {\'el\'ementaires} de la
forme
$$
\sum_{k_1\ge \cdots \ge k_p\ge 1} \frac{z_1^{-k_1}\cdots
z_p^{-k_p}} {(k_1+j_1)^{s_1}\cdots (k_p+j_p)^{s_p}},
$$
o\`u $s_i\in\mathbb{Z}$ et $j_i\in\mathbb{N}$.
C'est ce probl\`eme que nous  allons maintenant r\'esoudre ; cela terminera la preuve
du th\'eor\`eme \ref{theo:informel}.

\subsection{Notations}

Dans tout ce paragraphe,
$N$ d\'esignera un entier $\ge 1$ qui jouera essentiellement le r\^ole
de profondeur, r\^ole d\'evolu jusqu'\`a pr\'esent \`a l'entier $p$.
On notera :
\begin{itemize}
\item $\un{j}_N=(j_i)_{i=1,\ldots N}$ et $\un{m}_N=(m_i)_{i=1,\ldots N}$ (avec
$m_1=0$) des suites d'entiers de $\mathbb{N}$ ;
\item $\un{s}_N=(s_i)_{i=1,\ldots N}$ une suite d'entiers de $\mathbb{Z}$ ;
\item $\un{z}_N=(z_i)_{i=1,\ldots N}$ une suite de complexes de modules
$\ge 1$ ;
\item $a\wedge b=\min{(a,b)}$ et $a\vee b=\max{(a,b)}$ ;
\item $\varepsilon_{a,b}=1$ si $a<b$, $-1$ si
$a>b$, $0$ si $a=b$ et $\varepsilon_{p}=\varepsilon_{j_{p-1},j_p+m_p}$ (pour $p\ge 2$) ;
\item $t_p=j_{p-1}\wedge (j_p+m_p)$ et  $T_p=j_{p-1}\vee (j_p+m_p)$ (pour $p\ge 2$).
\end{itemize}
\medskip

\`A toute suite finie
$\un{u}_N=(u_1,\ldots,u_N)$, on associe les trois suites :
\begin{itemize}
\item $\un{u}_{N}^p=(u_p,\ldots,u_N)$ de longueur $N-p+1$ (pour $1\le p\le N$) ;
\item
$_p\un{u}_{N}=(u_1,\ldots,u_{p-2},u_{p-1}u_p,u_{p+1},\ldots,u_N)$ de longueur
$N-1$ (pour $2\le p\le N$) ;
\item $1/\un{u}_{N}=(1/u_1,\ldots,1/u_N)$ lorsque les $u_i$ sont non-nuls.
\end{itemize}

\medskip

On d\'efinit les {\em briques decal\'ees-modul\'ees} par
\begin{equation}
\ba_N\Bigg[\,
\begin{matrix}
\un{s}_N\\
\un{m}_N\\
\un{j}_N
\end{matrix}\,
\Bigg\vert \,\un{z}_N
\,\Bigg]=
\sum_{k_{N-1}+m_N \ge k_N \ge 1\atop
{k_{N-2}+m_{N-1} \ge k_{N-1} \ge 1\atop
{\vdots\atop
{k_1+m_2 \ge k_2\ge 1\atop
k_1\ge 1}}}
}\frac{z_1^{-k_1}\cdots z_N^{-k_N}}
{(k_1+j_1)^{s_1}\cdots(k_N+j_N)^{s_N}}.
\label{eq:3}
\end{equation}
Les $j_i$ sont les d\'ecalages, les
$m_i$ les modulations, les $s_i$ les exposants, $N$ la profondeur et on d\'efinit son poids comme \'etant $\sum_{p=1}^N \max(s_p,0)$.
Par d\'efinition,
$m_1=0$ : toutes les briques
$\ba^{\prime}$ que nous construirons \`a l'aide de briques $\ba$
avec $m_1=0$ auront aussi $m^{\prime}_1=0$.  Ces  s\'eries
convergent absolument lorsque $\vert z_1\vert>1$ et  $\vert z_j\vert \ge 1$ pour $j=2, \ldots, N$,
ce que l'on suppose
dor\'enavant et qui l\'egitime les diverses manipulations que nous
effectuerons dessus ; nous montrerons au paragraphe~\ref{sec:raffinement en z=1} comment obtenir des r\'esultats
similaires lorsque tous les $z_i$ valent 1.
Un cas particulier important est celui
o\`u tous les
$m_i$ sont nuls : on parlera de
{\em brique  d\'ecal\'ee}, ou simplement de {\em brique},~\footnote{Dans un contexte voisin, Zudilin~\cite{zud} a  introduit
une notion de {\em brique}, reprise et g\'en\'eralis\'ee dans~\cite{kratriv}. Ces briques n'ont rien \`a voir avec les n\^otres ; elles
sont suffisamment diff\'erentes pour ne pas les confondre si on est amen\'e \`a manipuler les deux types de briques
simultan\'ement.} et on la notera
\begin{equation}
\ba_N\bigg[\,
\begin{matrix}
\un{s}_N\\
\un{j}_N
\end{matrix}\,
\bigg\vert \,\un{z}_N\,
\bigg]=\sum_{k_1 \ge \cdots \ge k_{N}\ge 1}
\frac{z_1^{-k_1}\cdots z_N^{-k_N}}
{(k_1+j_1)^{s_1}\cdots(k_N+j_N)^{s_N}}.
\label{eq:4}
\end{equation}
Nous avons d\'ej\`a rencontr\'e ce type de briques dans les cas $N=1$ et $N=2$
au paragraphe~\ref{sec:prof 1 et 2} et en toute g\'en\'eralit\'e au
paragraphe~\ref{ssec:serie hyp mult en brique}.
Pour obtenir des relations compactes,
on d\'efinit la brique de profondeur 0 (et vide de param\`etres) comme la fonction identiquement \'egale \`a 1.
La modulation semble {\em a priori} une notion
artificielle et inutile puisqu'on ne s'int\'eresse r\'eellement qu'aux briques d\'ecal\'ees  :
\`a l'usage, il n'en est rien car, de fa\c{c}on surprenante, on ne
peut apparemment pas produire le th\'eor\`eme~\ref{prop:2} ci-dessous
sans modulation.

Nous appellerons {\em terme de profondeur $\le N-1$} toute
combinaison lin\'eaire \`a coefficients dans
$\mathbb{Q}[z^{\pm 1}_1,\ldots, z_N^{\pm 1}]$ de briques
d\'ecal\'ees-modul\'ees de profondeur $\le N-1$ et \'evalu\'ees
en des produits quelconques des variables $z_1, \ldots, z_N$. Le poids d'un terme de
profondeur $N-1$ est le plus grand des poids des briques qui le composent.

Pour tout entier $p$ tel que $1\le p\le N+1$, on d\'efinit le polyn\^ome de Laurent
\begin{equation*}
Q_{\un{s}_{N}^{p},p}(K;\un{z}_{N}^{p})
=\sum_{K\ge
k_{p}\ge \cdots\ge  k_{N}\ge 1}
\frac{z_{p}^{-k_{p}}\cdots z_{N}^{-k_{N}}}
{\prod_{i=p}^N k_i^{s_i}},
\label{eq:5}
\end{equation*}
(qui vaut $0$ si $K=0$)
pour $p\le N$ et $Q_{\un{s}_{N}^{N+1},N+1}(K;\un{z}_{N}^{N+1})=1$ pour $p=N+1$.
On notera
$Q_{\un{s}_{N}^p,p}(K;\un{z}_{N}^{p})=Q_{N,p}(K;\un{z}_{N}^{p})$
lorsqu'il n'y aura pas de risque de confusion sur les exposants en jeu.
On a
\begin{equation}
Q_{N,p}(K;\un{z}_{N}^{p})=
\sum_{k_{p}=1}^{K}\frac{z_p^{-k_p}}{k_p^{s_p}}\,
Q_{N,p+1}(k_{p};\un{z}_{N}^{p+1}).
\label{eq:6}
\end{equation}

Enfin, pour tout
entier $p$ tel que $2\le p\le N$, on d\'efinit
\begin{multline} R_{\un{s}_{N},p}(K; {}_p\un{z}_{N})
\\
=
\sum_{k_{p-2}+m_{p-1}\ge k_{p-1}\ge 1\atop{
\vdots\atop
{k_1+m_2\ge k_2\ge 1 \atop
k_1\ge 1}}}
\frac{z_1^{-k_1}\cdots z_{p-2}^{-k_{p-2}}(z_{p-1}z_p)^{-k_{p-1}}}
{\left(\prod_{i=1}^{p-1}(k_i+j_i)^{s_i}\right)
(k_{p-1}+K)^{s_p}} \,Q_{N,p+1}(k_{p-1}+K;\un{z}_{N}^{p+1}).
\label{eq:7}
\end{multline}
Si $p=2$, on attribue la valeur 1 au produit vide $z_1^{-k_1}\cdots z_{p-2}^{-k_{p-2}}$.
On notera $R_{\un{s}_{N},p}(K;{}_p\un{z}_{N})=R_{N,p}(K;{}_p\un{z}_{N})$ lorsqu'il
n'y aura pas de risque de confusion et nous
montrerons qu'il s'agit d'un terme de profondeur $\le
N-1$.

\subsection{L'algorithme de d\'ecomposition des briques}\label{sec:2}
Le but de ce paragraphe est de d\'emontrer que
{\em la brique d\'ecal\'ee-modul\'ee \eqref{eq:3} est la somme de
$\left(z_1^{j_1}\cdots z_N^{j_N}\right)\la_{\un{s}_N}(1/\un{z}_N)$
et de termes de profondeur au plus
$N-1$}. Cette proposition informelle (qui suffit \`a d\'emontrer le th\'eor\`eme \ref{theo:informel},
compte tenu des r\'esultats du paragraphe \ref{sec:decompositionenbrique}) d\'ecoule du th\'eor\`eme
suivant qui est beaucoup plus pr\'ecis.

\begin{theo} \label{prop:2}
$(i)$ Pour tout entier $N\ge 1$, on a
\begin{multline}
\ba_N\Bigg[\,
\begin{matrix}
\un{s}_N\\
\un{m}_N\\
\un{j}_N
\end{matrix}
\,\Bigg\vert\, \un{z}_N
\,\Bigg]=
\left(z_1^{j_1}\cdots
z_N^{j_N}\right)\la_{\un{s}_N}(1/\un{z}_N)
\\
-\sum_{p=1}^N \left(z_p^{j_p}\cdots
z_N^{j_N}\right)\,Q_{N,p}(j_p; \un{z}_{N}^{p})\,
\ba_{p-1}\Bigg[\,
\begin{matrix}
\un{s}_{p-1}\\
\un{m}_{p-1} \\
\un{j}_{p-1}
\end{matrix}\,
\Bigg\vert\,\un{z}_p\,
\Bigg]\\ \qquad+\sum_{p=2}^N \varepsilon_p\,
\left(z_p^{j_p}\cdots z_N^{j_N}\right)
\sum_{k_p=t_p +1}^{T_p}z_p^{-k_p}R_{N,p}(k_p; {}_p\un{z}_{N}).
\label{eq:8}
\end{multline}

$(ii)$ Pour tout entier $p$ tel que $2\le p\le N$ et tout entier $K\ge 0$, la s\'erie
$R_{N,p}(K;{}_p\un{z}_{N})$ est un terme de profondeur $\le N-1$, dont le poids est $\le \sum_{p=1}^N\max(s_p,0)$.
\end{theo}
\begin{remark} $(1)$ Si $N=1$, l'expression d\'ebutant par
$\sum_{p=2}^N\,\varepsilon_p (\cdots)$ n'appara{\^\i}t pas.

$(2)$ Ce th\'eor\`eme fournit
un algorithme permettant d'expliciter totalement le r\'esultat infor\-mel
 \'evoqu\'e au d\'ebut de ce paragraphe. Nous avons impl\'ement\'e cet algorithme sous
GP-Pari.
\end{remark}

\begin{proof}[D\'emonstration]

La partie (i) repose sur le lemme suivant, que nous
d\'emontrons \`a la toute fin
de ce paragraphe.

\begin{lem}\label{lem:2}
$(i)$ Pour tout $N\ge 2$ et tout $p=2,\ldots, N$,  on a
\begin{multline*}
\ba_N\Bigg[\,
\begin{array}{ccccc}
\un{s}_{N}& & & &\\
\un{m}_{p}& j_{p}& 0& \cdots&  0 \\
\un{j}_{p}& 0      & 0& \cdots&  0
\end{array}\,
\Bigg\vert \,\un{z}_N \,
\Bigg]
= z_p^{j_p}\,\ba_N\Bigg[\,
\begin{array}{ccccc}
\un{s}_{N}& & & &\\
\un{m}_{p-1}& j_{p-1}& 0& \cdots &  0 \\
\un{j}_{p-1}& 0      & 0& \cdots &  0
\end{array}\,
\Bigg\vert \,\un{z}_N\,
\Bigg]\\
-z_{p}^{j_p}\,Q_{N,p}(j_p; \un{z}_{N}^{p})\,
\ba_{p-1}\Bigg[\,
\begin{array}{c}
\un{s}_{p-1}\\
\un{m}_{p-1} \\
\un{j}_{p-1}
\end{array}\,
\Bigg\vert \,\un{z}_p\,
\Bigg]
+\varepsilon_p
\sum_{k_p=t_p +1}^{T_p} z_p^{j_p-k_p}R_{N,p}(k_p;{}_p\un{z}_{N}).
\end{multline*}

$(ii)$ Pour tout $N\ge 1$, on a
\begin{equation*}
\ba_N\Bigg[\,
\begin{array}{ccccc}
\un{s}_{N}& & & &\\
 0 & j_{1}& 0& \cdots&  0 \\
\un{j}_{1}& 0      & 0& \cdots&  0
\end{array}\,
\Bigg\vert\,\un{z}_N\,
\Bigg]=
z_1^{j_1}\,\la_{\un{s}_N}(1/\un{z}_N)-z_1^{j_1}\,Q_{N,1}(j_1;\un{z}_{N}).
\end{equation*}
\end{lem}

On applique le point (i) de ce lemme avec $p=N$, ce qui donne
\begin{multline*}
\ba_N\Bigg[\,
\begin{matrix}
\un{s}_N\\
\un{m}_N\\
\un{j}_N
\end{matrix}\,
\Bigg\vert \,\un{z}_N\,
\Bigg]=-z_N^{j_N} \,Q_{N,N}(j_N;z_N)\,
\ba_{N-1}\Bigg[\,
\begin{matrix}
\un{s}_{N-1}\\
\un{m}_{N-1}\\
\un{j}_{N-1}
\end{matrix}\,
\Bigg\vert \,\un{z}_{N-1}\,
\Bigg]\\
+z_N^{j_N}
\ba_N\Bigg[\,
\begin{array}{cc}
\un{s}_{N}& \\
\un{m}_{N-1}& j_{N-1}\\
\un{j}_{N-1}& 0
\end{array}\,
\Bigg\vert \,\un{z}_N\,
\Bigg ]
+\varepsilon_N
\sum_{k_N=t_N +1}^{T_N}z_N^{j_N-k_N}R_{N,N}(k_N; {}_N\un{z}_{N}).
\end{multline*}
On rep\`ete ce proc\'ed\'e $N-1$ fois en appliquant le
lemme \ref{lem:2}, (i),  \`a l'unique
brique de profondeur $N$ qui appara\^\i t \`a chaque it\'eration,
jusqu'\`a obtenir (en plus d'autres termes)
la brique
\begin{equation*}
\ba_N\Bigg[\,
\begin{array}{ccccc}
\un{s}_{N}&&&&\\
\un{m}_{1}& j_{1}& 0& \cdots & 0 \\
\un{j}_{1}& 0& 0&\cdots& 0
\end{array}\,
\Bigg\vert \,\un{z}_N
\,\Bigg],
\end{equation*}
\`a laquelle on applique alors le point (ii) du m\^eme
lemme \ref{lem:2} (puisque $\un{m}_{1}=m_1=0$).
En regroupant les termes, on
constate que l'on a d\'emontr\'e le point (i) du th\'eor\`eme~\ref{prop:2}.

\medskip

Pour prouver la partie (ii), on a \'egalement besoin d'un lemme technique, dont on
donnera la d\'emonstration \`a la fin du paragraphe.

\begin{lem}\label{lem:4} Soient $e, f \in\mathbb{Z}$ et $i, j\in\mathbb{C}$.

$(i)$ Lorsque $i=j$,
\begin{equation*}
\frac{1}{(X+i)^e(X+j)^f}=\frac{1}{(X+i)^{e+f}}\;.
\end{equation*}

$(ii)$ Lorsque $e\le 0$ et $f\ge 1$,
\begin{equation*}
\frac{1}{(X+i)^e(X+j)^f}=
\sum_{u=0}^{-e}\binom{-e}{u}(i-j)^{-e-u}\frac{1}{(X+j)^{f-u}}\;.
\end{equation*}

$(iii)$ Lorsque $e,f\le 0$,
\begin{equation*}
\frac{1}{(X+i)^e(X+j)^f}=
\sum_{u=0}^{-e}\sum_{v=0}^{-f}\binom{-e}{u}\binom{-f}{v}
i^{-e-u}j^{-f-v}X^{u+v}\;.
\end{equation*}

$(iv)$ Lorsque $i\not=j$ et $e,f\ge 1$,
\begin{equation*}
\frac{1}{(X+i)^e(X+j)^f}=\sum_{u=1}^e
\frac{\binom{e+f-1-u}{f-1}}{(i-j)^{e+f-u}}\frac{(-1)^f}{(X+i)^u}+
\sum_{v=1}^f\frac{\binom{e+f-1-v}{e-1}}{(j-i)^{e+f-v}}
\frac{(-1)^e}{(X+j)^v}.
\end{equation*}
\end{lem}

Nous allons exprimer $R_{N,p}(K;{}_p\un{z}_{N})$ en termes de briques
 \`a l'aide du lemme \ref{lem:4} appliqu\'e \`a la fraction
\begin{equation*}
\frac{1}{(k_{p-1}+j_{p-1})^{s_{p-1}}(k_{p-1}+K)^{s_p}}\;
\end{equation*}
qui appara\^\i t dans~\eqref{eq:7}.
Cinq
cas se pr\'esentent naturellement et
il n'y en a pas d'autres possibles ; leurs intersections peuvent
\^etre non vides mais c'est sans importance ici.

Si $p=N$, resp. $p=2$, les colonnes correspondant \`a
$s_{p+1}, s_{p+2}, \ldots, s_N$, resp. $\un{s}_{p-2}$, des six
briques $\ba_{N-1}$ suivantes
n'apparaissent pas.

\subsubsection{Premier cas : $K=j_{p-1}$} \label{sssec: 221}
Cela correspond au cas (i) du lemme \ref{lem:4}. On a alors
\begin{equation*}
R_{N,p}(K
;{}_p\un{z}_{N})
=\ba_{N-1}\Bigg[\,
\begin{array}{cccccc}
\un{s}_{p-2}&s_{p-1}+s_p&s_{p+1}&s_{p+2}&\cdots&s_N\\
\un{m}_{p-2}&m_{p-1}&j_{p-1}&0&\cdots&0\\
\un{j}_{p-2}&j_{p-1}&0&0&\cdots&0
\end{array}\,
\Bigg\vert \, {}_p\un{z}_{N}
\,\Bigg].
\end{equation*}

\subsubsection{Deuxi\`eme cas : $s_{p-1}\le 0$ et $s_p\ge 1$}
Cela correspond au cas (ii) du lemme \ref{lem:4}. On a alors
\begin{multline*}
R_{N,p}(K;{}_p\un{z}_{N})
=\sum_{u=0}^{-s_{p-1}}\binom{-s_{p-1}}{u} (j_{p-1}-K)^{-s_{p-1}-u}\\
\cdot\ba_{N-1}\Bigg[\,
\begin{array}{cccccc}
\un{s}_{p-2}&s_{p}-u &s_{p+1}&s_{p+2}&\cdots&s_N\\
\un{m}_{p-2}&m_{p-1}& K&0&\cdots&0\\
\un{j}_{p-2}& K& 0&0&\cdots&0
\end{array}\,
\Bigg\vert \, {}_p\un{z}_{N}\,
\Bigg].
\end{multline*}

\subsubsection{Troisi\`eme cas : $s_{p-1}\ge 1$ et $s_p\le 0$}
Cela correspond de nouveau au cas (ii) du lemme~\ref{lem:4}. On a alors
\begin{multline*}
R_{N,p}(K;{}_p\un{z}_{N})
=\sum_{u=0}^{-s_{p}}\binom{-s_{p}}{u}(K-j_{p-1})^{-s_{p}-u}\\
\cdot\ba_{N-1}\Bigg[\,
\begin{array}{cccccc}
\un{s}_{p-2}&s_{p-1}-u &s_{p+1}&s_{p+2}&\cdots&s_N\\
\un{m}_{p-2}&m_{p-1}& K&0&\cdots&0\\
\un{j}_{p-2}& j_{p-1}& 0&0&\cdots&0
\end{array}\,
\Bigg\vert \, {}_p\un{z}_{N}\,
\Bigg].
\end{multline*}

\subsubsection{Quatri\`eme cas : $s_{p-1}\le 0$ et $s_p\le 0$}
Cela correspond au cas (iii) du lemme~\ref{lem:4}. On a alors
\begin{multline*}
R_{N,p}(K;{}_p\un{z}_{N})
=\sum_{u=0}^{-s_{p-1}}\sum_{v=0}^{-s_{p}}
\binom{-s_{p-1}}{u}\binom{-s_{p}}{v}
j_{p-1}^{-s_{p-1}-u}K^{-s_{p}-v}\\
\cdot\ba_{N-1}\Bigg[\,
\begin{array}{cccccc}
\un{s}_{p-2}&-u-v&s_{p+1}&s_{p+2}&\cdots&s_N\\
\un{m}_{p-2}&m_{p-1}&K&0&\cdots&0\\
\un{j}_{p-2}&0&0&0&\cdots&0
\end{array}\,
\Bigg\vert \, {}_p\un{z}_{N}\,
\Bigg].
\end{multline*}

\subsubsection{Cinqui\`eme cas : $K\neq j_{p-1}$,  $s_{p-1}\ge 1$ et $s_p\ge 1$}
\label{sssec: 225}
Cela correspond au cas (iv) du lemme~\ref{lem:4}. On a alors
\begin{small}
\begin{multline*}
R_{N,p}(K;{}_p\un{z}_{N})\\
= (-1)^{s_{p}}
\sum_{u=1}^{s_{p-1}}
\frac{\binom{s_{p-1}+s_{p}-1-u}{s_{p}-1}}
{(j_{p-1}-K)^{s_{p-1}+s_{p}-u}}
\,\ba_{N-1}\Bigg[
\begin{array}{cccccc}
\un{s}_{p-2}&u&s_{p+1}&s_{p+2}&\cdots&s_N\\
\un{m}_{p-2}&m_{p-1}& K& 0&\cdots&0\\
\un{j}_{p-2}&j_{p-1}& 0&0&\cdots&0
\end{array}\,
\Bigg\vert \, {}_p\un{z}_{N}\,
\Bigg]
\\
+(-1)^{s_{p-1}}\sum_{v=1}^{s_p}
\frac{\binom{s_{p-1}+s_{p}-1-v}{s_{p-1}-1}}
{(K-j_{p-1})^{s_{p-1}+s_{p}-v}}
\,\ba_{N-1}\Bigg[\,
\begin{array}{cccccc}
\un{s}_{p-2}&v&s_{p+1}&s_{p+2}&\cdots&s_N\\
\un{m}_{p-2}&m_{p-1}& K& 0&\cdots&0\\
\un{j}_{p-2}&K& 0&0&\cdots&0
\end{array}\,
\Bigg\vert \, {}_p\un{z}_{N}\,
\Bigg].
\end{multline*}
\end{small}
Chacun de ces cinq cas montre que
$R_{N,p}(K;{}_p\un{z}_{N})$ est un terme de profondeur
$\le N-1$, de poids $\le \sum_{p=1}^N \max(s_p,0)$,
ce qui conclut la preuve.
\end{proof}

\begin{proof}[D\'emonstration du lemme~\ref{lem:2}]
Montrons (i). Remarquons tout d'abord que pour
toute suite $(u_n)_{n\ge 0}$, on a
\begin{eqnarray}
\lefteqn{\sum_{\ell=1}^{k+m}\frac{u_{\ell+j}\,z^{-\ell}}{(\ell+j)^s}
=\sum_{\ell=j+1}^{k+j+m}\frac{u_{\ell}\,z^{j-\ell}}{\ell^s}}\nonumber\\
&=&
\lm(-\sum_{\ell=1}^j+
\sum_{\ell=1}^{k+i}+\varepsilon_{i,j+m}
\sum_{\ell=k+i\wedge (j+m)+1}^{k+i\vee
(j+m)}
\rt)
\frac{u_{\ell}\,z^{j-\ell}}{\ell^s}\nonumber\\
&=&
\lm(-\sum_{\ell=1}^j+
\sum_{\ell=1}^{k+i}
\rt)
\frac{u_{\ell}\,z^{j-\ell}}{\ell^s}
+\varepsilon_{i,j+m}
\sum_{\ell=i\wedge (j+m)+1}^{i\vee
(j+m)}\frac{u_{k+\ell}\,z^{j-k-\ell}}{(k+\ell)^s},
\label{eq:suite}
\end{eqnarray}
apr\`es quelques manipulations imm\'ediates.

Supposons maintenant $2\le p\le N-1$.
On a
\begin{multline*}
\ba_N\Bigg[\,
\begin{array}{ccccc}
\un{s}_{N}&&&&\\
\un{m}_{p}&j_{p}&0&\cdots&0 \\
\un{j}_{p}& 0& 0&\cdots&0
\end{array}\,
\Bigg\vert \,\un{z}_N \,
\Bigg]=\\
\sum_{k_{p-2}+m_{p-1} \ge k_{p-1}\ge 1\atop{
\vdots\atop{
k_1+m_2 \ge k_2\ge 1\atop
k_1\ge 1}}} \frac{z_1^{-k_1}\cdots z_p^{-k_{p-1}}}
{\prod_{i=1}^{p-1}(k_i+j_i)^{s_i}}\sum_{k_p=1}^{k_{p-1}+m_p}
\frac{z_p^{-k_p}}{(k_p+j_p)^{s_p}}\,Q_{N,p+1}(k_p+j_p; \un{z}_{N}^{p+1})\,
\end{multline*}
On applique \eqref{eq:suite} \`a la somme $\sum_{k_p=1}^{k_{p-1}+m_p}(\ldots)$ et \`a
la suite
$u_{n}=Q_{N,p+1}(n; \un{x}_{N}^{p+1})$ : gr\^ace \`a la relation~\eqref{eq:6} entre
$Q_{N,p}$ et
$Q_{N,p+1}$, on voit alors que
\begin{multline*}
\ba_N\Bigg[\,
\begin{array}{ccccc}
\un{s}_{N}&&&&\\
\un{m}_{p}&j_{p}&0&\cdots&0 \\
\un{j}_{p}& 0& 0&\cdots&0
\end{array}\,
\Bigg\vert \,\un{z}_N\,
\Bigg]
= z_p^{j_p}\,\ba_N\Bigg[\,
\begin{array}{ccccc}
\un{s}_{N}&&&&\\
\un{m}_{p-1}&j_{p-1}&0&\cdots&0 \\
\un{j}_{p-1}& 0& 0&\cdots&0
\end{array}\,
\Bigg\vert \,\un{z}_N\,
\Bigg]\\
-z_{p}^{j_p}\,Q_{N,p}(j_p ; \un{z}_{N}^{p})\,
\ba_{p-1}\Bigg[\,
\begin{array}{c}
\un{s}_{p-1}\\
\un{m}_{p-1} \\
\un{j}_{p-1}
\end{array}
\Bigg\vert\,\un{z}_p\,
\Bigg]
+\varepsilon_p
\sum_{k_p=t_p +1}^{T_p}z_p^{j_p-k_p}R_{N,p}(k_p;{}_p\un{z}_{N}).
\end{multline*}

Pour (ii), on a
\begin{eqnarray*}
\lefteqn{
\ba_N\Bigg[\,
\begin{array}{ccccc}
\un{s}_{N}&&&&\\
0  &j_{1}&0&\cdots&0 \\
\un{j}_{1}& 0& 0&\cdots&0
\end{array}\,
\Bigg\vert \,\un{z}_N\,
\Bigg]
}\\
&=&\sum_{k_1=1}^{\infty}\frac{z_1^{-k_1}}
{(k_1+j_1)^{s_1}}\,Q_{N,2}(k_1+j_1;\un{z}_{N}^{2})
=\sum_{k_1=j_1+1}^{\infty}\frac{z_1^{j_1-k_1}}{k_1^{s_1}}\,Q_{N,2}(k_1;
\un{z}_{N}^{2})\\
&=& z_1^{j_1}\,\ba_N
\Bigg[\,
\begin{array}{ccccc}
\un{s}_{N}&&&&\\
0&0&0&\cdots&0 \\
0& 0& 0&\cdots&0
\end{array}\,
\Bigg\vert \,\un{z}_N\,
\Bigg]-\sum_{k_1=1}^{j_1}\frac{z_1^{j_1-k_1}}{k_1^{s_1}}\,
Q_{N,2}(k_1; \un{z}_{N}^{2})\\
&=& z_1^{j_1}\la_{\un{s}_N}(1/\un{z}_N)-z_1^{j_1}\,Q_{N,1}(j_1;\un{z}_{N})
\end{eqnarray*}
ce qui termine la d\'emonstration du lemme.
\end{proof}

\begin{proof}[D\'emonstration du lemme~\ref{lem:4}] Les points (i), (ii) et (iii) sont triviaux
et on d\'emontre seulement (iv), qui l'est \`a peine moins. En effet, on a
$$
\frac{1}{(X+i)^e(X+j)^f} = \sum_{u=1}^e \frac{a_u}{(X+i)^u} +
\sum_{v=1}^f \frac{b_v}{(X+j)^v}
$$
avec
\begin{multline*}
a_u = \frac{1}{(e-u)!} \bigg( \frac{1}{(X+j)^f}\bigg)^{e-u}\bigg\vert_{X=-i}
\\
=
\frac{(-f)(-f-1)\cdots (-f-e+u+1)}{(e-u)!(j-i)^{e+f-u}}=
\binom{e+f-u-1}{f-1}\frac{(-1)^f}{(i-j)^{e+f-u}}
\end{multline*}
et la formule similaire attendue pour $b_v.$
\end{proof}

\section{Pr\'ecisions sur le Th\'eor\`eme~\ref{prop:2}}\label{sec:raffinement algo}

Le but de ce paragraphe est de
pr\'eciser la nature des polyn\^omes de Laurent qui
apparaissent quand on it\`ere le Th\'eor\`eme \ref{prop:2}, sous la
condition que tous les exposants $s_i$ sont strictement positifs.

On pose
\begin{itemize}
\item $M_i=\sum_{k=1}^{i} m_k$ avec $M_0=0$ ;
\item $I_N=\max_{i=1, \ldots, N}(T_i+M_{i-1})$
avec $T_i=j_{i-1}\vee (j_i+m_i)$, $j_0=0$ et $I_0=0$ ;
\item $J_{N}=\max_{i=1,\ldots, N}(j_{i})$ et $J_0=0$ ;
\item $K_N=\max_{i=1, \ldots, N}(T_i)$ et $K_0=0$;
\item $\Sigma_N=\sum_{i=1}^N s_{i}$.
\end{itemize}
$J_N$ est le cas sp\'ecial de $I_N$ obtenu lorsque les modulations sont toutes nulles.
Rappelons que $\dd_n$ d\'enote le p.p.c.m. des entiers $1, 2, \ldots, n$. Par convention, $\dd_0=1$.
On utilisera le fait trivial que
$\dd_n^e \dd_m^f$ divise $\dd_{n\vee m}^{e+f}$.

\begin{theo}  \label{theo:raffinement} Supposons que tous
les exposants $s_i$ sont strictement positifs.

$(i)$ Les polyn\^omes de Laurent qui interviennent dans la d\'ecompo\-si\-tion de
la brique d\'ecal\'ee-modul\'ee large~\eqref{eq:3} en
polylogarithmes multiples
sont dans
\begin{equation*}
d_{I_{N}}^{-\Sigma_N}\mathbb{Z}
[z_1,z_2^{\pm 1},\ldots,z_N^{\pm 1}]
\end{equation*}
et leur degr\'e en $z_1$ est au plus $K_N$.

$(ii)$ Les polyn\^omes de Laurent qui interviennent dans la
d\'ecompo\-si\-tion de la brique d\'ecal\'ee large~\eqref{eq:4} en
polylogarithmes multiples sont dans
\begin{equation*}
d_{J_{N}}^{-\Sigma_{N}}\mathbb{Z}
\end{equation*}
et leur degr\'e en $z_1$ est au plus $J_N$.
\end{theo}

\begin{remark} $(1)$ Le point $(ii)$ est le seul vraiment utile ; nous ne
sa\-vons pas le d\'e\-mon\-trer sans d'abord d\'e\-mon\-trer $(i$), dont il est un cas particulier.

$(2)$ On n'utilisera pas que
les $s_i$ sont strictement positifs pour d\'emontrer
que les polyn\^omes de Laurent sont
des polyn\^omes de degr\'e $\le j_1$ en la variable $z_1$.

$(3)$ Concernant le d\'enominateur, un r\'esultat similaire  a probablement lieu dans le cas g\'en\'eral mais nous n'avons
pas cherch\'e \`a l'expliciter, faute de perspectives diophantiennes \'evidentes.
\end{remark}

\begin{proof} (i) Nous
proc\'edons, en deux temps, par r\'ecurrence sur la profondeur
$N$ de la brique
\eqref{eq:3} : le point (ii) en d\'ecoule en prenant le cas particulier de modulations toutes nulles.

\subsection{Preuve de l'assertion sur les d\'enominateurs}

Le cas $N=1$ est imm\'ediat : on a
\begin{equation}
\label{eq:B_1 recurrence}
\ba_1\Bigg[\,
\begin{array}{c}
s_{1}\\
0\\
j_{1}
\end{array}\,
\Bigg \vert \,z_1\,
\Bigg]=
z_1^{j_1}\,\la_{s_1}(z_1)-z_1^{j_1}\,Q_{1,1}(j_1;z_1),
\end{equation}
o\`u
$\displaystyle
Q_{1,1}(j_1;z_1)=\sum_{k_1=1}^{j_1}\frac{z_1^{-k_1}}
{k_1^{s_1}}$
a pour d\'enominateur $\dd_{j_1}^{s_1}=\dd_{I_1}^{\Sigma_1}$.

Supposons maintenant le Th\'eor\`eme~\ref{theo:raffinement} vrai
jusqu'\`a la profondeur $N-1$ et analysons les diff\'erents termes
de l'\'equation~\eqref{eq:8}, que nous rappelons :
\begin{multline*}
\ba_N\Bigg[\,
\begin{array}{c}
\un{s}_N\\
\un{m}_N\\
\un{j}_N
\end{array}\,
\Bigg\vert \,\un{z}_N\,
\Bigg]=
(z_1^{j_1}\cdots
z_N^{j_N})\,\la_{\un{s}_N}(\un{z}_N)\\
-\sum_{p=1}^N (z_p^{j_p}\cdots
z_N^{j_N}) \,Q_{N,p}(j_p; \un{z}_{N}^{p})\,\ba_{p-1}
\Bigg[\,
\begin{array}{c}
\un{s}_{p-1}\\
\un{m}_{p-1} \\
\un{j}_{p-1}
\end{array}\,
\Bigg\vert\, \un{z}_p\,
\Bigg]\\ +\sum_{p=2}^N \varepsilon_p\,(z_p^{j_p}\cdots z_N^{j_N})
\sum_{k_p=t_p +1}^{T_p}z_p^{-k_p}R_{N,p}(k_p; {}_p\un{z}_{N}).
\end{multline*}
Tout d'abord
\begin{equation*}
Q_{N,p}(j_p;\un{z}_{N}^p)
=\sum_{j_p\ge k_p\ge \cdots\ge k_N\ge 1}
\frac{z_{p}^{-k_{p}}\cdots z_{N}^{-k_{N}}}
{\prod_{i=p}^N k_i^{s_i}}
\end{equation*}
a pour d\'enominateur $\dd_{j_p}^{s_p+\cdots+s_N}$. Par hypoth\`ese de
r\'ecurrence, un d\'enomina\-teur de la brique
$\ba_{p-1}$ est $\dd_{I_{p-1}}^{s_1+\cdots+s_{p-1}}$, m\^eme pour $p=1$. Un
d\'enominateur des termes
$Q_{N,p}\ba_{p-1}$ est donc
$\dd_{j_p}^{s_p+\cdots+s_N}\,\dd_{I_{p-1}}^{s_1+\cdots+s_{p-1}}$,
qui divise $\dd_{I_{N}}^{\Sigma_N}$ puisque $j_p\vee
I_{p-1}\le (T_p+M_{p-1})\vee
I_{p-1} = I_p\le  I_N$ pour tout $p\in\{1,\ldots, N\}$.

Il reste \`a analyser les termes $R_{N,p}(k_p; {}_p\un{z}_{N})$ : nous allons distinguer
deux cas.

\subsubsection{Premier cas : $k_p=j_{p-1}$}
On est alors dans la situation du paragraphe~\ref{sssec: 221}:
\begin{equation*}
R_{N,p}(j_{p-1};{}_p\un{z}_{N})
=\ba_{N-1}\Bigg[\,
\begin{array}{cccccc}
\un{s}_{p-2}&s_{p-1}+s_p&s_{p+1}&s_{p+2}&\cdots&s_N\\
\un{m}_{p-2}&m_{p-1}&j_{p-1}&0&\cdots&0\\
\un{j}_{p-2}&j_{p-1}&0&0&\cdots&0
\end{array}\,
\Bigg\vert \, {}_p\un{z}_{N}\,
\Bigg].\label{eq:10}
\end{equation*}
L'hypoth\`ese de r\'ecurrence s'applique : un d\'enominateur de la
brique  est
$$
\dd_{I_{p-1}\vee (j_{p-1}\vee (0+j_{p-1})+M_{p} )}^{\Sigma_N} =
\dd_{I_{p-1}\vee (j_{p-1} +M_{p} )}^{\Sigma_N}.
$$
Comme
$j_{p-1}+M_{p-1}\le T_{p}+M_{p-1}$, ce d\'enominateur divise $\dd_{I_{p}}^{\Sigma_N}$,
qui divise $\dd_{I_{N}}^{\Sigma_N}$.

\subsubsection{Second cas : $k_p\not=j_{p-1}$}
On est maintenant dans la
situation du paragraphe~\ref{sssec: 225} :
\begin{equation*}
\sum_{p=2}^N\varepsilon_p\,(z_p^{j_p}\cdots z_N^{j_N})\sum_{k_p=t_p +1\atop
k_p\not=j_{p-1}}^{T_p}R_{N,p}(k_p;{}_p\un{z}_{N})
=\sum_{p=2}^N\varepsilon_p\,(z_p^{j_p}\cdots z_N^{j_N})\bigg(\sum_{u=1}^{s_{p-1}}B_{1,p}(u)
+\sum_{v=1}^{s_p}B_{2,p}(v)\bigg)
\end{equation*}
avec
\begin{multline}
B_{1,p}(u)=(-1)^{s_{p}}\sum_{k_p=t_p +1\atop k_p\not=j_{p-1}}^{T_p}
\frac{\binom{s_{p-1}+s_{p}-1-v}{s_{p}-1}}{(j_{p-1}-k_p)^{s_{p-1}+s_{p}-u}}\\
\cdot
\ba_{N-1}\Bigg[\,
\begin{array}{cccccc}
\un{s}_{p-2}&u&s_{p+1}&s_{p+2}&\cdots&s_N\\
\un{m}_{p-2}&m_{p-1}& k_p& 0&\cdots&0\\
\un{j}_{p-2}&j_{p-1}& 0&0&\cdots&0
\end{array}\,
\Bigg\vert\, {}_p\un{z}_{N}
\Bigg]  \label{eq:cpu}
\end{multline}
et
\begin{multline}
B_{2,p}(v)=(-1)^{s_{p-1}}\sum_{k_p=t_p
+1\atop
k_p\not=j_{p-1}}^{T_p}
\frac{\binom{s_{p-1}+s_{p}-1-v}{s_{p-1}-1}}
{(k_p-j_{p-1})^{s_{p-1}+s_{p}-v}}\\
\cdot \ba_{N-1}\Bigg[\,
\begin{array}{cccccc}
\un{s}_{p-2}&v&s_{p+1}&s_{p+2}&\cdots&s_N\\
\un{m}_{p-2}&m_{p-1}& k_p& 0&\cdots&0\\
\un{j}_{p-2}&k_p& 0&0&\cdots&0
\end{array}\,
\Bigg\vert\, {}_p\un{z}_{N}\,
\Bigg].
\label{eq:dpu}
\end{multline}
Nous allons montrer que $\dd_{I_{N}}^{\Sigma_N}$ est un d\'enominateur
convenable pour les termes \eqref{eq:cpu} et
\eqref{eq:dpu}, ce qui suffira puisqu'il est ind\'ependant de $p$, $u$ et $v$.
Fixons
$p$,
$u$ et
$v$. Par hypoth\`ese de r\'ecurrence, les deux briques
$\ba_{N-1}$ ont pour  d\'enominateurs respectifs
$$
D_1 = \dd_{I_{p-1}\vee (k_p+ M_{p-1})}^{u+\Sigma_N-s_p-s_{p-1}}
\quad \textup{et} \quad
D_2 = \dd_{I_{p-2}\vee (j_{p-2}\vee(k_p+m_{p-1})+M_{p-2})\vee
(k_p +M_{p-1})}^{v+\Sigma_N-s_p-s_{p-1}}.
$$
Puisque $k_p\le T_p$, on a
$I_{p-1}\vee (k_p+M_{p-1})\le I_{p-1}\vee (T_p + M_{p-1}) = I_p$ et donc $D_1$ divise
$\dd_{I_p}^{u+\Sigma_N-s_p-s_{p-1}}$. D'autre part, si $j_{p-2}\le k_p+ m_{p-1}$, on a
$$
j_{p-2}\vee(k_p+m_{p-1})+M_{p-2}\le k_p+m_{p-1}+M_{p-2}\le T_p+M_{p-1}
$$ tandis que si
$j_{p-2}\ge k_p+ m_{p-1}$,
alors
$$
j_{p-2}\vee(k_p+m_{p-1})+M_{p-2}\le j_{p-2}+M_{p-2}\le T_{p-1}+M_{p-2},
$$
d'o\`u $D_2$ divise  $\dd_{I_{p-2}\vee (T_{p-1}+M_{p-2})\vee
(T_p +M_{p-1})}^{v+\Sigma_N-s_p-s_{p-1}}= \dd_{I_p}^{v+\Sigma_N-s_p-s_{p-1}}$.
On obtient donc des
d\'enomina\-teurs uniformes en
$k_p$ pour les briques $\ba_{N-1}$ :
$$
\dd_{I_{p}}^{u+\Sigma_N-s_p-s_{p-1}}
\quad \textup{et}\quad
\dd_{I_{p}}^{v+\Sigma_N-s_p-s_{p-1}}.
$$
Les deux sommes
$$
\dd_{I_{p}}^{u+\Sigma_N-s_p-s_{p-1}}B_{1,p}(u)=(-1)^{s_p}\sum_{k_p=t_p +1\atop k_p\not=j_{p-1}}^{T_p}
\frac{\binom{s_{p-1}+s_{p}-1-u}{s_{p}-1}}
{(j_{p-1}-k_p)^{s_{p-1}+s_{p}-u}}\,\dd_{I_{p}}^{u+
\Sigma_N-s_p-s_{p-1}}\ba_{N-1}[\cdots]
$$
et
$$
\dd_{I_{p}}^{v+\Sigma_N-s_p-s_{p-1}} B_{2,p}(v)=(-1)^{s_{p-1}}\sum_{k_p=t_p
+1\atop
k_p\not=j_{p-1}}^{T_p}
\frac{\binom{s_{p-1}+s_{p}-1-v}{s_{p-1}-1}}
{(k_p-j_{p-1})^{s_{p-1}+s_{p}-v}}\,
\dd_{I_{p}}^{v+\Sigma_N-s_p-s_{p-1}}
\ba_{N-1}[\cdots]
$$
ont donc pour d\'enominateurs respectifs
$
\dd_{\vert j_p+m_p-j_{p-1}\vert}^{-u+s_{p-1}+s_p}
$
et
$
\dd_{\vert j_p+m_p-j_{p-1}\vert}^{-v+s_{p-1}+s_p},
$
qui divisent trivialement
$
\dd_{I_{p}}^{-u+s_{p-1}+s_p}$, resp. $\dd_{I_{p}}^{-v+s_{p-1}+s_p}
$,
car $\vert j_p+m_p-j_{p-1}\vert\le T_p\le I_p$.
Ainsi, on peut prendre
\begin{equation*}
\dd_{I_{p}}^{-u+s_{p-1}+s_p}\,
\dd_{I_{p}}^{u+\Sigma_N-s_p-s_{p-1}}=\dd_{I_{p}}^{\Sigma_N}
\label{eq:den1}
\end{equation*}
et
\begin{equation*}
\dd_{I_{p}}^{-v+s_{p-1}+s_p}
\,\dd_{I_{p}}^{v+\Sigma_N-s_p-s_{p-1}}
=\dd_{I_{p}}^{\Sigma_N},
\label{eq:den2}
\end{equation*}
comme d\'enominateur de \eqref{eq:cpu} et \eqref{eq:dpu}, ce qui ach\`eve la
preuve du  Th\'eor\`eme \ref{theo:raffinement} puisque $\dd_{I_p}^{\Sigma_N}$ divise
$\dd_{I_N}^{\Sigma_N}$.

\subsection{Preuve de l'assertion sur le degr\'e en $z_1$}

De nouveau, on raisonne par r\'ecurrence sur la profondeur $N\ge 1$.
C'est \'evidemment vrai pour $N=1$
par l'\'equation~\eqref{eq:B_1 recurrence}.
Supposons maintenant l'assertion vraie pour $N-1$ et, comme pr\'ec\'edemment, analysons les
termes de l'\'equation~\eqref{eq:8}.
Le terme $(z_1^{j_1}\cdots
z_N^{j_N})\,\la_{\un{s}_N}(\un{z}_N)$ est de la forme voulue, avec un degr\'e
$j_1\le K_N$.
Dans le terme
$$
\sum_{p=1}^N (z_p^{j_p}\cdots
z_N^{j_N}) \,Q_{N,p}(j_p; \un{z}_{N}^{p})\,\ba_{p-1}
\Bigg[\,
\begin{array}{c}
\un{s}_{p-1}\\
\un{m}_{p-1} \\
\un{j}_{p-1}
\end{array}\,
\Bigg\vert\, \un{z}_p\,
\Bigg],
$$
si $p\ge 2$, la variable $z_1$ n'appara{\^\i}t pas dans les polyn\^omes de
Laurent  $Q_{N,p}(j_p; \un{z}_{N}^{p})$ et seulement dans la brique $B_{p-1}[\ldots]$ qui est
de profondeur $p-1\le N-1$ : l'hypoth\`ese de r\'ecurrence s'applique et seules les puissances positives de
$z_1$ interviennent bien, jusqu'au plus $z_1^{K_{p-1}}$, donc au plus
$z_1^{K_N}$. Si $p=1$, alors $z_1$ intervient dans l'expression
$z_1^{j_1}\,Q_{N,1}(j_1; \un{z}_{N}^{1})\,\ba_{0}[\ldots] = z_1^{j_1}\,Q_{N,1}(j_1; \un{z}_{N})$,
qui est aussi un polyn\^ome en $z_1$ de degr\'e au plus $j_1 \le K_N$.
Il reste le dernier terme
$$
\sum_{p=2}^N \varepsilon_p\,(z_p^{j_p}\cdots z_N^{j_N})
\sum_{k_p=t_p +1}^{T_p}z_p^{-k_p}R_{N,p}(k_p; {}_p\un{z}_{N})
$$
qui ne d\'epend de $z_1$ que par $R_{N,p}(k_p; {}_p\un{z}_{N})$. Or les expressions que nous
en avons donn\'ees au paragraphe pr\'ec\'edent montrent
qu'il s'agit d'une combinaison lin\'eaire de briques
de profondeur $\le N-1$ \'evalu\'ees en ${}_p\un{z}_{N}$ et dont les coefficients ne d\'ependent pas des
$z_i$. Dans ${}_p\un{z}_{N}$, la variable
$z_1$ appara{\^\i}t seule si $3 \le p\le N$ : l'hypoth\`ese de r\'ecurrence s'applique
et on v\'erifie que le degr\'e en $z_1$ est au plus
$K_{p}\le K_N$. Si $p=2$, alors il y a une subtilit\'e car  $z_1$ appara{\^\i}t
multipli\'e par $z_2$ : ce n'est pas g\^enant, l'hypoth\`ese de r\'ecurrence s'applique
de nouveau et le degr\'e en $z_1$ est $\le T_2\le
K_N$,
ce qui conclut la d\'emonstration.
\end{proof}

\section{Non-enrichissement des $\la_{s_1,\ldots, s_p}$ \`a exposants n\'egatifs} \label{sec:nonenrichissement}

L'algorithme de d\'ecomposition des briques peut faire appara\^itre des polylogarithmes larges \`a exposants n\'egatifs (ou nuls). Par
exemple la d\'ecomposition de l'int\'egrale de Sorokin pour $\zeta (3)$
$$
\int_{[0,1]^3} \di \frac{u^n (1-u)^n v^n (1-v)^n w^n (1-w)^n}{(z_1 -uv)^{n+1}
(z_1z_2 -uvw)^{n+1}} \,\dd u \dd v \dd w ,
$$
fait intervenir des $\la_{s_1,s_2} (1/z_1, 1/z_2 )$, avec $s_1=1,2$, $s_2=0, -1, \dots, -n+1$.

Afin de r\'egler ces cas singuliers, on d\'emontre un r\'esultat dit de {\it non-enrichissement arithm\'etique}.
\begin{theo}
\label{enri} Supposons que, pour tout $j=1, \ldots, p$, on ait $\vert z_j\vert <1$. Alors,
tout polyloga\-ri\-thme multiple large $\la_{s_1,\dots ,s_p} (z_1 ,\dots ,z_p)$
de profondeur $p$ ayant certains exposants $s_j\le 0$ s'exprime comme une combinaison lin\'eaire finie
de polylogarithmes multi\-ples lar\-ges $\la_{s'_1,\dots ,s'_q} (z_1^* ,\dots ,z_q^* )$
de profondeur $q \in \{0, \dots, p\}$, avec $s'_j\ge 1$,
o\`u les $z_i^*$
sont certains produits des $z_j$.
Les coefficients de la combinaison lin\'eaire sont des
polyn\^omes \`a coefficients rationnels
en les $\di \big((1-z_{j_1}\cdots z_{j_m})^{-1}\big)_{1\le j_1< \ldots < j_m\le p, \,m\ge 1}$ et les
$\big(z_j^{\pm 1}\big)_{1\le j\le p}$. De plus, on a $\sum_{j=1}^q s'_j\le \sum_{j=1}^p \max(0,s_j)$ pour toutes les
suites d'exposants $\underline{s}'$ qui apparaissent.
\end{theo}

\begin{remark} $(1)$ Pour tout $z$ tel que $\vert z\vert <1$, on a
$$
\la_{-s} (z) = \left(z\frac{\dd}{ \dd z}\right)^s \left(\frac{1}{1-z}\right) \in (1-z)^{-s-1}\mathbb{Z}[z].
$$

$(2)$ Ce th\'eor\`eme est de facture informelle mais sa d\'emonstration
offre un moyen algorith\-mique de l'expliciter.

$(3)$ Un r\'esultat de ce type est annonc\'e par \'Ecalle (\cite[pp. 419--420]{ecalle}) dans le cas des poly\-z\^etas, sans
d\'emonstration.
\end{remark}

\subsection{Pr\'eliminaires}

On suppose dans toute la suite de ce paragraphe que toutes les variables not\'ees $z$ ou $z_j$ sont de modules $<1$.
La d\'emonstration utilisera l'identit\'e triviale suivante, valable pour tout entier $K\ge 1$ :
\begin{equation}
\sum_{k_1 =1}^{K} \sum_{k_2=1}^{k_1} = \sum_{k_2 =1}^{K}
\bigg( \sum_{k_1=1}^{K} -\sum_{k_1=1}^{k_2 -1} \bigg ).
\label{eq:star}
\end{equation}
Pour tous entiers $s\ge 0$ et $K\ge 1$, on d\'efinit
$P_{s} (K,z)=\di\sum_{k=1}^K k^s z^k,$ qui v\'erifie :
\begin{equation*}\label{eq:pszN}
P_{s} (K,z) =
\bigg(z\frac{\dd}{\dd z}\bigg)^s \left (z \frac{1-z^{K}}{1-z} \right ).
\end{equation*}
On en d\'eduit que l'on a
\begin{equation}\label{eq:pszNzneq1}
P_{s} (K,z)= \sum_{\ell=0}^s \frac{z^Ka_{1, \ell}(s,z)+a_{2,\ell}(s,z)}{(1-z)^{s+1}}
K^\ell
\end{equation}
o\`u  $a_{1, \ell}(s,z)$ et $a_{1, \ell}(s,z)$ sont des polyn\^omes en $z$
de degr\'e au plus $s$ et ind\'ependants de $K$.   On notera

Les objets naturels qui vont intervenir sont des {\it polylogarithmes larges tronqu\'es} :
$$
\la^K_{s_1 , \dots , s_p}
(z_1 ,\dots ,z_p )= \di\sum_{K\geq k_1 \geq \dots \geq
k_p \geq 1} \di \frac{z_1^{k_1} \dots z_p^{k_p}}{k_1^{s_1}
\dots k_p^{s_p} }.
$$
On remarque que l'on a $\la^K_{s_1} (z_1) =  P_{-s_1}(K,z_1)$ lorsque $s_1\le 0$.

On aura besoin du lemme suivant.

\begin{lem}
\label{utile}
Soient des entiers $s_1\ge 0$ et $s_2, \ldots ,s_p\in\mathbb{Z}$. Pour tous entiers
$K\geq 1$ et $p\geq 2$, on a  :
\begin{multline*}
\la^K_{-s_1 ,s_2 , \dots , s_p } (z_1 ,\dots ,z_p ) =
P_{s_1}(K,z_1)
\la^K_{s_2 , \dots , s_p } (z_2 ,\dots ,z_p )
\\
- \sum_{\ell=0}^{s_1} \frac{a_{1,\ell}(s_1, z_1)}{(1-z_1)^{s_1+1}z_1}\sum_{m=0}^\ell
\binom{\ell}{m}
(-1)^{\ell-m} \la^K_{s_2 -m , s_3, \dots , s_p} (z_1z_2 , z_3, \ldots ,z_p )
\\
- \sum_{\ell=0}^{s_1} \frac{a_{2,\ell}(s_1, z_1)}{(1-z_1)^{s_1+1}}\sum_{m=0}^\ell
\binom{\ell}{m}
(-1)^{\ell-m} \la^K_{s_2 -m , s_3, \dots , s_p} (z_2 ,\dots ,z_p ).
\end{multline*}
\end{lem}

\begin{proof}
En utilisant \eqref{eq:star}, on a :
\begin{eqnarray*}
\la^K_{-s_1 ,s_2 , \dots , s_p } (z_1 ,\dots ,z_p )
&=&\sum_{k_2 =1}^K \frac{z_2^{k_2}}{k_2^{s_2}}
\bigg( \sum_{k_1 =1}^K k_1^{s_1} z^{k_1} - \di\sum_{k_1
=1}^{k_2 -1} k_1^{s_1} z^{k_1} \bigg)
\la^{k_2 }_{s_3  , \dots , s_p } (z_3 ,\dots ,z_p )
\\
&=&
\sum_{k_2 =1}^K {z_2^{k_2} \over k_2^{s_2}} \big(
P_{s_1} (K,z_1) - P_{s_1} (k_2 -1,z_1) \big)
\la^{k_2}_{s_3 , \dots , s_p } (z_3 ,\dots ,z_p ) .
\end{eqnarray*}
Au moyen de~\eqref{eq:pszNzneq1}, on obtient
\begin{multline*}
\la^K_{-s_1 ,s_2 , \dots , s_p } (z_1 ,\dots ,z_p )
=P_{s_1}(K,z_1)
\sum_{k_2=1}^K \bigg({z_2^{k_2} \over k_2^{s_2}}
\la^{k_2}_{s_3 , \dots , s_p} (z_3 ,\dots ,z_p ) \bigg)
\\
-
\sum_{\ell=0}^{s_1}
\frac{1}{(1-z_1)^{s_1+1}}
\sum_{k_2=1}^K (z_1^{k_2-1}a_{1, \ell}(s_1,z_1)+a_{2,\ell}(s_1,z_1))\bigg({z_2^{k_2} \over k_2^{s_2}}
\,(k_2 - 1)^\ell\,
 \la^{k_2}_{s_3 , \dots , s_p} (z_3 ,\dots ,z_p ) \bigg).
\end{multline*}
La premi\`ere somme vaut exactement
$$
P_{s_1}(K,z_1)
\la^K_{s_2 , \dots , s_p} (z_2 ,\dots ,z_p ).
$$
La seconde somme faisant intervenir $(k_2 -1)^\ell$ est \`a peine plus
compliqu\'ee. En d\'eveloppant le terme $(k_2 -1)^\ell$ par le th\'eor\`eme binomial et en rempla\c{c}ant
directement dans la somme, on obtient en effet :
\begin{multline*}
- \sum_{\ell=0}^{s_1} \frac{a_{1,\ell}(s_1, z_1)}{(1-z_1)^{s_1+1}z_1}\sum_{m=0}^\ell
\binom{\ell}{m}
(-1)^{\ell-m} \la^K_{s_2 -m , s_3, \dots , s_p} (z_1z_2 , z_3, \dots ,z_p )
\\
- \sum_{\ell=0}^{s_1} \frac{a_{2,\ell}(s_1, z_1)}{(1-z_1)^{s_1+1}}\sum_{m=0}^\ell
\binom{\ell}{m}
(-1)^{\ell-m} \la^K_{s_2 -m , s_3,\dots , s_p} (z_2 , z_3, \dots ,z_p ),
\end{multline*}
ce qui termine la d\'emonstration.
\end{proof}

\subsection{D\'emonstration du th\'eor\`eme~\ref{enri}}

On remarque que le lemme~\ref{utile} exprime un
polylogarithme de profondeur $p$ \`a l'aide de polylogarithmes
de profondeur $p-1$, ce qui ouvre la porte \`a une d\'emonstration du  th\'eor\`eme~\ref{enri} par r\'ecurrence

Pour $p=1$, le th\'eor\`eme est vrai, comme le montre la remarque (1) qui suit son \'enonc\'e.

On suppose que l'on sait d\'ecomposer les polylogarithmes de profondeur $\leq p-1$ (avec $p-1\ge 1$)
de la mani\`ere pr\'evu par le
th\'eor\`eme.
Soit maintenant $s_1 ,\dots , s_p$  une suite quelconque d'entiers, avec au moins un $s_j\le 0$ :
notons ${q+1}$ le plus petit indice $\ge 1$ tel que $s_{q+1} \le 0$.
Pour simplifier, on note $s_{q+1} =-s$ avec $s \ge 0$. On doit distinguer trois cas : $q=0$,  $1\le q\le p-2$ et $q=p-1.$

-- Le cas $q=0$. Notons que pour tout entier $t\ge 0$, on a
$$
\sum_{k=\ell}^{\infty} k^{t} z^{k} =
\left(z\frac{\dd}{ \dd z}\right)^t \left(\frac{z^{\ell}}{1-z}\right) = \frac{z^\ell Q_{t}(\ell,z)}{(1-z)^{t+1}}
$$
avec $Q_{t}(\ell,z)\in\mathbb{Z}[\ell,z]$ de degr\'e $s$ en $\ell$ et $z$. On pose donc
$
Q_t(\ell,z) = \sum_{j=0}^s q_{j,s}(z) \ell^j.
$
On a alors
\begin{eqnarray*}
\la_{s_1 ,s_2 , \dots , s_p } (z_1 , z_2, \dots ,z_p )
&=&  \sum_{k_2 \ge \cdots \ge k_p\ge 1} \bigg(
\frac{z_2^{k_2}\cdots z_p^{k_p}}{k_2^{s_2}\cdots k_p^{s_p}} \sum_{k_1=k_2}^{\infty} k_1^{s} z_1^{k_1}\bigg)\\
&=&\frac{1}{(1-z_1)^{s+1}} \sum_{k_2 \ge \cdots \ge k_p\ge 1} Q_{s}(k_2, z_1)
\frac{(z_1z_2)^{k_2}z_3^{k_3}\cdots z_p^{k_p}}{k_2^{s_2}k_3^{s_3}\cdots k_p^{s_p}}\\
&=& \frac{1}{(1-z_1)^{s+1}} \sum_{j=0}^s q_{j,s}(z_1) \la_{s_2 -j, s_3, \dots , s_p } (z_1z_2, z_3 ,\dots ,z_p).
\end{eqnarray*}
Comme on n'a finalement que des $\la$ de profondeur $p-1$, l'hypoth\`ese de r\'ecurrence s'ap\-pli\-que.

\medskip

-- Le cas $1\le q\le p-2$. On applique le lemme~\ref{utile} de telle sorte que
\begin{eqnarray}
\lefteqn{\la_{s_1 ,s_2 , \dots , s_p } (z_1 , z_2, \dots ,z_p ) }\nonumber
\\
&=&
\sum_{k_1\ge \cdots\ge k_q\ge 1}
\frac{z_1^{k_1}\cdots z_q^{k_q}}{k_1^{s_1}\cdots k_q^{s_q}} \,
\la^{k_q}_{-s, s_{q+2}, \ldots, s_p} (z_{q+1}, z_{q+2}, \ldots, z_p)\nonumber
\\
&=& \sum_{k_1\ge \cdots \ge k_q\ge 1}
\frac{z_1^{k_1}\cdots z_q^{k_q}}{k_1^{s_1}\cdots k_q^{s_q}} P_{s}(k_q,z_{q+1})
\la^{k_q}_{s_{q+2} , \dots , s_p } (z_{q+2} ,\dots ,z_p )\label{eq:tanguy1}
\\
&& \quad -
\sum_{\ell=0}^{s} \frac{a_{1,\ell} (s ,z_{q+1} )}{(1-z_{q+1})^{s+1}z_{q+1}} \sum_{m=0}^\ell \bigg((-1)^{\ell-m} \binom{\ell}{m}
\nonumber \\
&& \qquad  \qquad\times  \sum_{k_1\ge \cdots \ge k_q\ge 1}
\frac{z_1^{k_1}\cdots z_q^{k_q}}{k_1^{s_1}\cdots k_q^{s_q}} \,
\la^{k_q}_{s_{q+2} - m , s_{q+3} , \dots , s_p } (z_{q+1}z_{q+2} , z_{q+3} ,\dots ,z_p )\bigg) \label{eq:tanguy22}\\
&& \qquad
-\sum_{\ell=0}^{s} \frac{a_{2,\ell} (s ,z_{q+1} )}{(1-z_{q+1})^{s+1}} \sum_{m=0}^\ell \bigg((-1)^{\ell-m} \binom{\ell}{m}
\nonumber \\
&& \qquad \qquad \times  \sum_{k_1\ge \cdots \ge k_q\ge 1}
\frac{z_1^{k_1}\cdots z_q^{k_q}}{k_1^{s_1}\cdots k_q^{s_q}} \,
\la^{k_q}_{s_{q+2} - m , s_{q+3} , \dots , s_p } (z_{q+2} , z_{q+3} ,\dots ,z_p )\bigg)
\label{eq:tanguy2}
\end{eqnarray}

Il est facile de traiter les s\'eries~\eqref{eq:tanguy22} et~\eqref{eq:tanguy2} puisqu'elles valent respectivment
$$
\la_{s_1, \ldots, s_{q}, s_{q+2}-m, s_{q+3}, \ldots, s_p} (z_1, \ldots, z_{q}, z_{q+1}z_{q+2}, \ldots, z_p)
$$
et
$$
\la_{s_1, \ldots, s_{q}, s_{q+2}-m, s_{q+3}, \ldots, s_p} (z_1, \ldots, z_{q}, z_{q+2}, \ldots, z_p),
$$
qui sont de profondeur $p-1$ : on peut donc leur appliquer l'hypoth\`ese de r\'ecurrence.

Reste la s\'erie sur la ligne~\eqref{eq:tanguy1} : on utilise de nouveau la forme
d\'evelopp\'ee~\eqref{eq:pszNzneq1} de $P_s(K,z)$ pour en obtenir l'expression alternative
\begin{multline*}
\sum_{\ell=0}^s \sum_{k_1\ge \cdots \ge k_q\ge 1}
\frac{z_1^{k_1}\cdots z_q^{k_q}}{k_1^{s_1}\cdots k_q^{s_q}}
\frac{z_{q+1}^{k_q} a_{1, \ell}(s,z_{q+1})+a_{2,\ell}(s,z_{q+1})}{k_q^{-\ell}(1-z_{q+1})^{s+1}} \,
\la^{k_q}_{s_{q+2} , \dots , s_p } (z_{q+2} ,\dots ,z_p )
\\
=\frac{1}{(1-z_{q+1})^{s+1}}\sum_{\ell=0}^s \bigg(a_{1, \ell} (s,z_{q+1}) \la_{s_1, \ldots, s_{q-1},
s_{q}-\ell, s_{q+2}, \ldots, s_p}(z_1, \ldots, z_{q}z_{q+1}, z_{q+2}, \ldots, z_p)
\\
+ a_{2, \ell} (s,z_{q+1}) \la_{s_1, \ldots, s_{q-1},
s_{q}-\ell, s_{q+2}, \ldots, s_p}(z_1, \ldots, z_{q}, z_{q+2}, \ldots, z_p) \bigg).
\end{multline*}
Comme  on a maintenant affaire \`a
une combinaison lin\'eaire de $\la$ de profondeur $p-1$, l'hypoth\`ese de r\'ecurrence s'ap\-pli\-que.

\medskip

-- Le cas $q =p-1$. On a
\begin{eqnarray*}
\lefteqn{\la_{s_1 ,s_2 , \dots , s_p } (z_1 , z_2, \dots ,z_p ) }\nonumber
\\
&=&
\sum_{k_1\ge \cdots\ge k_{p-1}\ge 1}
\frac{z_1^{k_1}\cdots z_q^{k_{p-1}}}{k_1^{s_1}\cdots k_{p-1}^{s_{p-1}}} \,
\la^{k_{p-1}}_{-s} (z_p)\nonumber
\\
&=& \sum_{k_1\ge \cdots \ge k_{p-1}\ge 1}
\frac{z_1^{k_1}\cdots z_q^{k_q}}{k_1^{s_1}\cdots k_q^{s_q}} P_{s}(k_{p-1},z_{p})
\\
&=& \sum_{\ell=0}^s \sum_{k_1\ge \cdots \ge k_{p-1}\ge 1}
\frac{z_1^{k_1}\cdots z_q^{k_q}}{k_1^{s_1}\cdots k_q^{s_q}}
\frac{z_{p}^{k_{p-1}} a_{1, \ell}(s,z_{p})+a_{2,\ell}(s,z_{p})}{ k_q^{-\ell}(1-z_{p})^{s+1}}
\\
&=& \frac{1}{(1-z_{p})^{s+1}} \sum_{\ell=0}^s
\bigg( a_{1, \ell}(s,z_{p})\la_{s_1, \ldots, s_{p-2}, s_{p-2}-\ell}(z_1, \ldots, z_{p-2}, z_{p-1}z_p)
\\
&& \qquad +
a_{1, \ell}(s,z_{p})\la_{s_1, \ldots, s_{p-2}, s_{p-2}-\ell}(z_1, \ldots, z_{p-2}, z_{p-1})
\bigg).
\end{eqnarray*}
On peut de nouveau appliquer l'hypoth\`ese  de r\'ecurrence, ce qui termine la preuve du th\'e\-o\-r\`e\-me~\ref{enri}.

\section{D\'emonstration du th\'eor\`eme \ref{theocv}}

\label{sec:raffinement en z=1}

Pour d\'emontrer le th\'eor\`eme \ref{theocv}, nous devons r\'egulariser les
polyz\^etas divergents intervenant dans la d\'ecomposition d'une brique. La r\'egularisation  qui s'impose ici est
la r\'egularisation dite {\it shuffle} des polyz\^etas bas\'ee sur l'\'etude du comportement asymptotique des polylogarithmes
lorsque $z$ tend vers $1$.

\subsection{R\'egularisation $\shu$ analytique}

Dans \cite[Corollaire 2.5]{ra}, Racinet caract\'erise, suivant
les travaux de L. Boutet de Monvel, le comportement asymptotique
des polylogarithmes lorsque $z$ tend vers $1$.
\begin{theo}
\label{theo:asymp} Pour tous entiers strictement positifs $s_1
,\dots ,s_p$, la fonction $\li_{s_1 ,\dots ,s_p} (z)$ admet,
lorsque $z$ tend vers $1$ tel que $\vert z\vert <1$, un d\'eveloppement
asymptotique du type
$$\li_{s_1 ,\dots ,s_p} (z) =Q_{s_1 ,\dots ,s_p} ( \log (1-z) )+o ( (1-z)^{\varepsilon} )$$
avec $Q_{s_1 ,\dots ,s_p} \in \mathbb{C} [t ]$ et $\varepsilon \in
\mathbb{R}^*_+$.
\end{theo}
On note $\zeta^{\shu} (s_1 ,\dots ,s_p )$ la valeur
r\'egularis\'ee de $\zeta (s_1 ,\dots ,s_p )$ pour $s_1 =1$
obtenue en posant $\zeta^{\shu} (s_1 ,\dots ,s_p ) =Q(0)$, i.e. le
terme constant du polyn\^ome $Q$. Si $s_1 \geq 2$, on a bien s\^ur
$\zeta^{\shu} (s_1 ,\dots ,s_p ) =\zeta (s_1 ,\dots ,s_p )$.

\subsection{Aspects effectifs}

Notons que l'impl\'ementation effective de l'algorithme de d\'ecomposition demande deux choses:

\medskip

(i) Le calcul des $\zeta^{\shu} (s_1 ,\dots ,s_p )$ r\'egularis\'es en fonction des $\zeta$ classiques.

\medskip

(ii) Le calcul explicite du reste intervenant dans l'estimation asymptotique du
th\'eo\-r\`eme~\ref{theo:asymp}.

\medskip

\noindent Le calcul des $\zeta^{\shu} (s_1 ,\dots ,s_p )$ dans le
cas divergents peut s'effectuer de fa\c{c}on {\it combinatoire},
beaucoup plus simple que {\em via} le calcul effectif des
d\'eveloppements asymptotiques du th\'eo\-r\`eme~\ref{theo:asymp}.

\subsection{R\'egularisation $\shu$ combinatoire}

La r\'egularisation $\shu$ que nous venons de d\'efinir conserve
la sym\'etrie $\shu$ v\'erifi\'ee par les polyz\^etas convergents.
Soit $A=\{ \mathbf{0},\mathbf{1} \}$ un alphabet. On note $A_c$
l'ensemble des mots de $A^*$ commencant par $\mathbf{0}$ et se
terminant par $\mathbf{1}$. On note $\pi$ le morphisme de
$\mathbb{R} \langle A_c \rangle $ dans $\mathbb{R} \langle Y
\rangle$ d\'efini par
$\pi (\mathbf{0}^{s-1} \mathbf{1} )=y_s,$
pour tout $s\geq 1$. On note encore $\zeta$ le morphisme d\'efini
sur $A_c$ par $\zeta (\mathbf{0}^{s-1} \mathbf{1} )=
\zeta_s$.
Le produit de battage ou shuffle sur $A$ se d\'efinit par
r\'ecurrence sur la longueur des mots par
$$a\mathbf{b}\shu c\mathbf{d} =a (\mathbf{b}\shu c\mathbf{d}) +c(a\mathbf{b}\shu \mathbf{d}) ,$$
pour tout mot $\mathbf{b},\mathbf{d} \in A^*$, $a,c\in A$.

On d\'emontre en utilisant l'\'ecriture int\'egrale des
polyz\^etas la relation de sym\'etrie dite {\it shuffle} : Pour
tout $\mathbf{u} \in A_c$, $\mathbf{v}\in A_c$, on a
\begin{equation}
\label{symetrie:shu} \zeta (\mathbf{u} )\zeta (\mathbf{v} )=\zeta
(\mathbf{u} \shu \mathbf{v} )
\end{equation}
On renvoie par exemple \`a l'article~\cite{colmez} pour plus de d\'etails.

On note $A_0$ l'ensemble des mots de $A^*$ se terminant par
$\mathbf{1}$. On peut donner un sens aux polyzetas sur $A_0$ en
utilisant la relation (\ref{symetrie:shu}) en supposant que
celle-ci est encore v\'erifi\'ee pour tout mot de $A_0$, ce qui
est le cas de la r\'egularisation $\zeta^{\shu}$ ci-dessus. On
note encore $\zeta^{\shu}$ le polyzeta
\'etendu \`a $A_0$.

Pour tout mot $\bs =s_1 \dots s_r \in A_c$, $r\geq 1$, $s_i \in
A$, on a
$\mathbf{1} \shu \mathbf{1}^{i} \bs =(i+1)\mathbf{1}^{i+1} \bs +\mathbf{1}^i s_1 [\mathbf{1} \shu \bs^{>1} ],$
o\`u $\bs^{>1} =s_2 \dots s_r$. En appliquant $\zeta^{\shu}$, on
obtient
\begin{equation}
\label{recurshuffle} \zeta^{\shu} (\mathbf{1} ) \zeta^{\shu}
(\mathbf{1}^i \bs )=(i+1)\zeta^{\shu} (\mathbf{1}^{i+1} \bs ) +
\zeta^{\shu} (\mathbf{1}^i s_1 [\mathbf{1}\shu \bs^{>1} ]).
\end{equation}
Il est donc possible de calculer $\zeta^{\shu} (\mathbf{1}^{i+1}
\bs )$ par r\'ecurrence sur le nombre de $\mathbf{1}$.
Pour cela, il suffit de fixer une valeur \`a $\zeta^{\shu} (\mathbf{1} )$.

Pour obtenir une r\'egularisation combinatoire qui coincide avec
la r\'egularisation analytique d\'efinie au paragraphe
pr\'ec\'edent, on doit poser $\zeta^{\shu} (1) =0$. En effet, un
simple calcul donne $\li_1 (z)=-\log (1-z)$. La formule
(\ref{recurshuffle}) permet alors le calcul explicite et
algorithmique des polyz\^etas divergents.

\subsection{\'Enonc\'es}

Dans cette partie, et dans toute la suite, on pose pour   $j \in \{1, \ldots, p\}$ :
$$
D_j = \Big( \sum_{i=1} ^j A_i (n_i+1) \Big) - j - 1.
$$
\begin{lem} \label{lemCNSCV}
La s\'erie
$$\sum_{k_1\ge \cdots \ge k_p\ge 1}
\frac{P(k_1,  \ldots, k_p)}{(k_1)_{n_1+1}^{A_1}
 \cdots (k_p)_{n_p+1}^{A_p}}$$
 converge si, et seulement si, le polyn\^ome $P(X_1, \ldots, X_p)$ v\'erifie
 \begin{equation} \label{eqCV}
\sum_{i=1}^j \deg_{X_i}P \leq D_j \mbox{ pour tout } j \in \{1,
\ldots, p\}.
\end{equation}
 \end{lem}

\begin{remark} Lorsque $n_1 = \dots = n_p =0$ et $P=1$, ce lemme donne
les conditions exactes de convergence des polyz\^etas $\zeta(A_1, A_2, \ldots, A_p)$
lorsque les $A_j$ sont dans $\mathbb{Z}.$ Elles correspondent bien aux conditions
qui assurent la convergence absolue des polyz\^etas pour des exposants {\em complexes}.
Voir~\cite[p. 10]{kratriv2} pour une preuve de ces conditions.
\end{remark}

\begin{proof}
Pour d\'emontrer ce lemme, on va montrer en fait que  les
conditions \eqref{eqCV} \'equivalent au fait que, pour tout $B
\geq 0$, la s\'erie
\begin{equation} \label{eqCVsergal}
\sum_{k_1\ge \cdots \ge k_p\ge 1} \frac{P(k_1,  \ldots, k_p) (\log
k_p)^B}{(k_1)_{n_1+1}^{A_1}
\cdots (k_p)_{n_p+1}^{A_p}}
\end{equation}
converge. C'est \'evident pour $p=1$, puisque
les conditions  \eqref{eqCV} se r\'eduisent alors \`a $\deg_{X_1}P \leq A_1 (n_1+1)-2$. Supposons que
ce soit vrai pour $p-1$, et soit $P(X_1, \ldots, X_p)$ ;
posons $\delta = \deg_{X_p} P$. Si $\delta \leq A_p (n_p+1)-1$ alors on a
$$1 \ll \sum_{k_p = 1}^{k_{p-1}} \frac{k_p ^{\delta}
(\log k_p)^B}{ (k_p)_{n_p+1}^{A_p}} \ll (\log k_{p-1})^{B+\delta}$$
donc la convergence de \eqref{eqCVsergal} \'equivaut \`a celle de
\eqref{eqCVsergal} en profondeur $p-1$. Comme justement l'\'equation
correspondant \`a $j=p$ dans \eqref{eqCV} se d\'eduit des autres
(puisqu'on a suppos\'e $\deg_{X_p} P \leq A_p (n_p+1)-1$), la preuve
est termin\'ee dans ce cas. Supposons maintenant que l'on ait
$\delta \geq A_p (n_p+1)$. Alors on a
$$k_{p-1}^{\delta - A_p (n_p+1) + 1}   \ll \sum_{k_p = 1}^{k_{p-1}}
\frac{k_p ^{\delta} (\log k_p)^B}{ (k_p)_{n_p+1}^{A_p}}
\ll k_{p-1}^{\delta - A_p (n_p+1) + 1}   (\log k_{p-1})^{B}$$
donc la convergence de \eqref{eqCVsergal} avec $P(X_1, \ldots, X_p)$
\'equivaut \`a celle de  \eqref{eqCVsergal} avec un polyn\^ome
$\Pti(X_1, \ldots, X_{p-1})$ v\'erifiant $\deg_{X_i}\Pti = \deg_{X_i} P$
pour $i \in \{1, \ldots, p-2\}$ et $\deg_{X_{p-1}}\Pti = \deg_{X_{p-1}}
P + \deg_{X_p}P - A_p (n_p+1) + 1$. Or justement les conditions~\eqref{eqCV} pour un tel
 polyn\^ome $\Pti$ \'equivalent aux conditions~\eqref{eqCV} pour $P$.
Le lemme est donc d\'emontr\'e.
\end{proof}

On dit qu'une fonction $f$, d\'efinie sur un ouvert dont le point
1 appartient \`a l'adh\'erence, est {\em \`a divergence au plus
logarithmique} en $z=1$ si elle admet un d\'eveloppement
asymptotique de la forme
$f(z) = Q( \log(1-z)) + \gdo((1-z)^\eps)$
pour un certain $\eps > 0$ et un polyn\^ome $Q \in \C[t]$. La {\em
valeur r\'egularis\'ee} de $f$ en 1 est le coefficient constant de
$Q$, c'est-\`a-dire $Q(0)$. Dans le cas particulier o\`u $f$ est
d\'efinie et continue en 1, le polyn\^ome $Q$ est constant et
cette valeur r\'egularis\'ee est simplement $f(1)$.
 \begin{lem} \label{lemfctlog}
 La fonction
 \begin{equation} \label{eqfctlog}
\sum_{k_1\ge \cdots \ge k_p\ge 1} \frac{P(k_1,  \ldots,
k_p)}{(k_1)_{n_1+1}^{A_1}
 \cdots (k_p)_{n_p+1}^{A_p}} z^{-k_1}
 \end{equation}
est \`a divergence au plus logarithmique en $z=1$
 si, et seulement si,
 \begin{equation} \label{eqdivlog}
\sum_{i=1}^j \deg_{X_i}P \leq D_j + 1 \mbox{ pour tout } j \in
\{1, \ldots, p\}.
\end{equation}
 \end{lem}
La preuve de ce lemme est analogue \`a celle du lemme
\ref{lemCNSCV} ; seule l'initialisation diff\`ere vraiment,
puisque la fonction $\sum_{k \geq 1} k^{-1} z^{-k}$ \`a est
divergence au plus logarithmique en $z=1$.

Pour d\'emontrer le th\'eor\`eme \ref{theocv}, on va en fait d\'emontrer le r\'esultat suivant qui est plus fort.

\begin{theo} \label{theoreg}
Si les relations \eqref{eqdivlog} sont satisfaites alors la valeur
r\'egularis\'ee en 1 de la fonction \eqref{eqfctlog}  est   une
combinaison lin\'eaire \`a coefficients rationnels en les
polyz\^etas r\'egularis\'es $\zeta^{\shu} (s_1 ,\dots ,s_q )$ o\`u
$0\leq q\leq p$, $s_i \geq 1$ pour  $i=1,\dots ,q$, $\sum_{j=1}^q s_j \leq \sum_{j=1}^p A_j$.
En outre, on peut calculer {\em explicitement} une telle
combinaison lin\'eaire.
\end{theo}

\subsection{Preuve du th\'eor\`eme~\ref{theoreg}} \label{ssec:preuvetheoreg}

On d\'emontre le th\'eor\`eme \ref{theoreg} par r\'ecurrence sur
la profondeur $p$. Quand $p=0$, ce th\'eor\`eme est trivial ; les
arguments qui suivent permettent de le d\'emontrer pour $p=1$,
mais un raisonnement direct est beaucoup plus facile dans ce cas.
Supposons donc que ce th\'eor\`eme soit vrai en toute profondeur
strictement inf\'erieure \`a $p$.

Soit $P(X_1, \ldots, X_p)$ un polyn\^ome  tel que les relations
\eqref{eqdivlog} soient satisfaites. On pose
$$
R(X_1, \ldots, X_p)= \frac{P(X_1,  \ldots,
X_p)}{(X_1)_{n_1+1}^{A_1}
 \cdots (X_p)_{n_p+1}^{A_p}},
$$
et on \'etudie la fonction
$$f(z) = \sum_{k_1\ge \cdots \ge k_p\ge 1} R(k_1,  \ldots, k_p) z^{-k_1}$$
qui est d\'efinie pour $|z|  > 1$ et est  \`a divergence au plus
logarithmique en $z=1$ gr\^ace au lemme \ref{lemfctlog}.  On
utilise le d\'eveloppement en \'el\'ements simples de $R$, comme
au paragraphe \ref{ssec:serie hyp mult en brique} (dont on
reprend les notations). Ceci permet d'\'ecrire, pour $|z| > 1$ :
\begin{equation} \label{eqsomtot}
f(z) = \sum_{\qqq} C[\qqq] \sum_{k_1\ge \cdots \ge k_p\ge 1}
\frac{\prod_{i \in I} k_i^{\sich}}{\prod_{i \in
\Ic}(k_i+j_i)^{s_i}}z^{-k_1}.
\end{equation}
Dans cette formule et dans toute la suite, on note $\qqq$ un
quadruplet g\'en\'erique
$$(I, (s_i)_{i \in \Ic}, (j_i)_{i \in \Ic}, (\sich)_{i \in I} )$$
tel que $1 \leq s_i \leq A_i $ et $0 \leq j_i \leq n_i$ pour tout
$i \in \Ic$, et $0 \leq \sich \leq \aich $ pour tout $i \in I$. On
pose alors $C[\qqq] = \Ciandco$.

La difficult\'e est que ce d\'eveloppement en \'el\'ements simples
fait appara\^{\i}tre des fonctions de $z$ dont la divergence en 1
n'est pas logarithmique. Par exemple, si $p=2$, $n_1 = 2$, $A_1 =
1$, $P(X_1, X_2) = (X_2)_{n_2+1}^{A_2} X_2$ alors les  relations
\eqref{eqdivlog} sont satisfaites mais dans l'expression
\eqref{eqsomtot} apparaissent les sommes
$$\sum_{k_1 \geq k_2 \geq 1} \frac{k_2}{k_1+j} z^{-k_1}$$
pour $j \in \{0,1,2\}$, qui sont chacune \`a divergence non
logarithmique. Une m\'ethode pour r\'esoudre ce probl\`eme serait
de g\'en\'eraliser le th\'eor\`eme \ref{theoreg}, en autorisant
des divergences non logarithmiques (c'est-\`a-dire des
d\'eveloppements asymptotiques avec des termes
$\frac{\log^k(z)}{(1-z)^\ell}$). Mais cela n\'ecessiterait   une
g\'en\'eralisation du th\'eor\`eme \ref{theo:asymp}, et ne
pr\'esenterait pas d'int\'er\^et pratique. En effet, la pr\'esence
de p\^oles en $\frac{1}{1-z}$ n\'ecessite de conna\^{\i}tre aussi
le coefficient de $1-z$ dans les d\'eveloppements asymptotiques,
car leur produit contribue \`a la valeur en $z=1$. L'algorithme
devrait donc calculer beaucoup de termes des d\'eveloppements
asymptotiques, ce qui serait co\^uteux en temps et en m\'emoire.
C'est pourquoi on proc\`ede plut\^ot comme suit. L'id\'ee importante
est celle de la r\'egularisation : quand seules des divergences
logarithmiques  sont pr\'esentes,  seul le coefficient constant du
polyn\^ome en $\log(1-z)$ intervient dans les calculs, y compris
lorsqu'on doit faire des produits.

Notons $\cale_0$ l'ensemble des quadruplets $\qqq  = (I, (s_i)_{i
\in \Ic}, (j_i)_{i \in \Ic}, (\sich)_{i \in I} )$ tels que la
fonction
\begin{equation} \label{equnebriq}
 \sum_{k_1\ge \cdots \ge k_p\ge 1} \frac{\prod_{i \in I} k_i^{\sich}}{\prod_{i \in \Ic}(k_i+j_i)^{s_i}}z^{-k_1}
\end{equation}
soit  \`a divergence au plus logarithmique en $z=1$,  et $\cale_1$
son compl\'ementaire. Dans la somme \eqref{eqsomtot}, chaque
\'el\'ement $\qqq \in \cale_0$ donne lieu \`a un d\'eveloppement
asymptotique de la forme
$Q_{\qqq}(\log(1-z)) + \gdo((1-z)^\eps)$
avec $\eps > 0$ (qu'on peut choisir ind\'ependant de $\qqq$) et
$Q_{\qqq} \in \C[t]$. En regroupant d'autre part les contributions
de tous les \'el\'ements $\qqq \in \cale_1$, on a donc :
\begin{equation} \label{eq1235}
 \sum_{\qqq \in \cale_1} C[\qqq] \sum_{k_1\ge \cdots \ge k_p\ge 1}
\frac{\prod_{i \in I} k_i^{\sich}}{\prod_{i \in \Ic}
(k_i+j_i)^{s_i}}z^{-k_1} = f(z) -  \sum_{\qqq \in \cale_0} Q_{\qqq}(\log(1-z)) + \gdo((1-z)^\eps).
\end{equation}
Comme $f(z)$ est \`a divergence au plus  logarithmique, on voit
que le membre de gauche aussi ; on va maintenant transformer ce
membre de gauche en une somme du type \eqref{eqfctlog} en
profondeur $p-1$. Soit $\qqq \in \cale_1$, avec $\qqq = (I,
(s_i)_{i \in \Ic}, (j_i)_{i \in \Ic}, (\sich)_{i \in I} )$.
L'hypoth\`ese \eqref{eqdivlog} (avec $j=1$) montre que l'ensemble
$J$ d\'efini au paragraphe \ref{ssec:serie hyp mult en brique} est
inclus dans $\{2, \ldots, p \}$, donc $I$ aussi. En outre, $I$ est
non vide (sinon on aurait $\qqq \in \cale_0$ d'apr\`es le lemme
\ref{lemfctlog}). Donc il  existe $t \in \{2, \ldots, p\}$ tel que
$t \in I$. En notant $B_s$ le $s$-i\`eme polyn\^ome de Bernoulli
(qui est \`a coefficients rationnels), on a \footnote{On peut noter que l'on utilise les m\^emes id\'ees
que celles du paragraphe~\ref{sec:nonenrichissement} sur le non-enrichissement des $\la$ \`a exposants n\'egatifs.
En particulier, \eqref{eq1237} est l'analogue de~\eqref{eq:pszNzneq1} lorsque tous les $z_j$ valent 1.} :
\begin{equation} \label{eq1237}
\sum_{k_t = k_{t+1} } ^{k_{t-1}} k_t^{\stch} = B_{\stch}
(k_{t-1}+1) - B_{\stch}(k_{t+1}).
\end{equation}
Cette relation permet d'\'ecrire, en posant $\ell_1 = k_1$,
\ldots, $\ell_{t-1} = k_{t-1}$, $\ell_t = k_{t+1}$, $\ell_{p-1} =
k_p$ :
\begin{multline*}
\sum_{k_1\ge \cdots \ge k_p\ge 1} \frac{\prod_{i \in I} k_i^{\sich}}{\prod_{i \in \Ic}(k_i+j_i)^{s_i}}z^{-k_1}\\
=\sum_{\ell_1\ge \cdots \ge \ell_{p-1} \ge 1} \frac{
\displaystyle{\prod_\indso{i \in I}{i \leq t-1} \ell_i^{\sich}
\prod_\indso{i \in I}{i \geq t+1} \ell_{i-1}^{\sich} } }{
\displaystyle{\prod_\indso{i \in \Ic}{i \leq
t-1}(\ell_i+j_i)^{s_i} \prod_\indso{i \in \Ic}{i \geq
t+1}(\ell_{i-1}+j_i)^{s_i}}
 } \times \Big(  B_{\stch} (\ell_{t-1}+1) - B_{\stch}(\ell_t) \Big)  z^{-\ell_1} .
 \end{multline*}
Cette somme est de la forme
$$\sum_{\ell_1\ge \cdots \ge \ell_{p-1} \ge 1} R_{\qqq} (\ell_1, \ldots, \ell_{p-1}) z^{-\ell_1}$$
pour une certaine fraction rationnelle $R_{\qqq}$ (qui d\'epend
aussi du choix, arbitraire et fix\'e, de~$t$). Le membre de gauche
de \eqref{eq1235} s'\'ecrit donc
\begin{equation} \label{eq1236}
\sum_{\ell_1\ge \cdots \ge \ell_{p-1} \ge 1} \Rtilde(\ell_1,
\ldots, \ell_{p-1}) z^{-\ell_1},
\end{equation}
o\`u l'on a pos\'e
$$ \Rtilde(\ell_1, \ldots, \ell_{p-1}) =  \sum_{\qqq \in \cale_1} C[\qqq]  R_{\qqq} (\ell_1, \ldots, \ell_{p-1}).$$
La relation \eqref{eq1235} et le lemme \ref{lemfctlog} montrent
que cette fraction rationnelle $ \Rtilde(\ell_1, \ldots,
\ell_{p-1})$ satisfait aux hypoth\`eses du th\'eor\`eme
\ref{theoreg}, en profondeur $p-1$.  Par hypoth\`ese de
r\'ecurrence, on peut donc \'ecrire \eqref{eq1236} sous la forme
$\Qtilde (\log(1-z)) + \gdo((1-z)^\eps)$, o\`u $\Qtilde(0)$ est
une   combinaison lin\'eaire explicite \`a coefficients rationnels
en les polyz\^etas r\'egularis\'es $\zeta^{\shu} (s_1 ,\dots ,s_q
)$ o\`u $1\leq q\leq p-1$, $\sum_{j=1}^q s_j \leq \sum_{j=1}^p A_j$. Compte tenu de \eqref{eq1235}, il suffit
maintenant de calculer $Q_{\qqq}(0)$ pour $\qqq \in \cale_0$, et
la preuve du th\'eor\`eme \ref{theoreg} sera termin\'ee.

Pour cela, on d\'ecompose la somme \eqref{equnebriq}. Tout
d'abord, si $I$ est non vide alors on  applique la relation
\eqref{eq1237} comme ci-dessus, et on est ramen\'e \`a une
profondeur strictement inf\'erieure. On peut donc supposer que $I$
est vide. Il suffit alors de suivre la preuve du th\'eor\`eme
\ref{theo:informel} (voir le paragraphe \ref{sec:l'algorithme}) avec $z_1 = z$, $z_2 =
\ldots = z_p = 1$, puis d'appliquer le
th\'eor\`eme \ref{theoreg} en profondeur $\leq p-1$. Ceci termine
la preuve du th\'eor\`eme \ref{theoreg}.


\def\refname{Bibliographie}

\bigskip

J. Cresson,
Laboratoire de Math\'ematiques appliqu\'ees de Pau,
B\^atiment I.P.R.A, Univer\-si\-t\'e de Pau et des Pays de l'Adour,
avenue de l'Universit\'e, BP 1155, 64013 Pau cedex, France.

\medskip

S. Fischler,
\'Equipe d'Arithm\'etique et de G\'eom\'etrie Alg\'ebrique,
Universit\'e Paris-Sud,
B\^atiment 425,
91405 Orsay Cedex, France.

\medskip

T. Rivoal,
Institut Fourier,
CNRS UMR 5582, Universit{\'e} Grenoble 1,
100 rue des Maths, BP~74,
38402 Saint-Martin d'H{\`e}res cedex,
France.

\end{document}